\documentclass[11pt]{article}

\usepackage{microtype}  
\usepackage{sectsty}      
\allsectionsfont{\sffamily}
\usepackage[margin=15pt,font=small,labelfont={bf,sf}]{caption}  

\newcommand{\s}[1]{{\huge\textsf{\textbf{#1}}}}

\usepackage{amsmath}
\usepackage{amsfonts}
\usepackage[pdftex]{graphicx}
\usepackage{breqn}
\usepackage{float}
\usepackage{pdfpages}
\usepackage{multirow}
\usepackage{longtable}
\usepackage{hyperref}
\usepackage[capitalise,sort&compress]{cleveref}
\hypersetup{colorlinks=true,urlcolor=black, linkcolor=black,citecolor=blue}

\usepackage[utf8]{inputenc}
\usepackage[T1]{fontenc}
\usepackage{textgreek}

\usepackage{mathtools}

\usepackage{epstopdf}
\usepackage{manfnt}
\reversemarginpar
\usepackage{xcolor}
\usepackage{bm}
\usepackage{caption}
\usepackage{subcaption}
\usepackage{wrapfig} 
\usepackage{textcomp}

\graphicspath{{./figures/}}

\title{\s{A network aggregation model for the dynamics and treatment of neurodegenerative diseases at the brain scale} }
\author{
Georgia S. Brennan{$^*$}
 and Alain Goriely{$^*$}\\ \\
$^*$\small{Mathematical Institute, Andrew Wiles Building}\\
\small{Woodstock Rd, University of 
Oxford, OX2 6GG, UK}
}

\begin{document}
\maketitle
\begin{abstract}\vskip 6pt
\hrule\vskip 12pt
Neurodegenerative diseases are associated with  the assembly of specific proteins into oligomers and fibrillar aggregates. At the brain scale, these protein assemblies can diffuse through the brain and seed other regions, creating an autocatalytic protein progression.
The growth and transport of these assemblies depend on various mechanisms that can be targeted therapeutically. Here, we use spatially-extended nucleation-aggregation-fragmentation models for the dynamics of prion-like neurodegenerative protein-spreading in the brain to study the effect of different drugs on whole-brain Alzheimer's disease progression.
\vskip 12pt
\hrule
\end{abstract}

\section{Introduction}

The self-assembly of proteins into ordered linear
structures plays a central role in the normal functioning of organisms, spanning from bacteria to mammals \cite{fowler2006functional, maji2009functional}, with unique mechanical properties also eliciting potential in industrial contexts \cite{bleem2017structural,knowles2016amyloid}. Recent interest in protein aggregation has led to an explosion in exploratory studies across a broad spectrum of disciplines. A particularly pressing motive for such research has more medical origins;
the undesired filamentous aggregation of proteins can have severe repercussions on an organism's well-being. In such contexts, aggregation is the defining mechanism driving a cascade of pathogenic proteins characteristic of various diseases. There are now approximately 50 disorders associated with a particular class of protein filaments, known as amyloids, with disparate symptoms ranging from non-neuropathic localized amyloidosis, such as type II diabetes, to neurodegenerative diseases such as Alzheimer's disease (AD), Parkinson's disease and Huntington's disease \cite{chiti2006protein}. Despite the diversity of pathogenic proteins involved in these disorders \cite{chiti2017protein}, experimental studies uncover commonalities in the underlying physicochemical and biochemical disease origins. This shared characteristic is the misfolding of normally soluble, functional peptides and proteins and their subsequent conversion into intractable aggregates \cite{knowles2014amyloid}. 
Although only understood just decades ago neurodegenerative diseases are no longer rare, and are rapidly becoming among the most common and debilitating medical conditions in the modern world \cite{
wortmann2012dementia}. This growing problem poses significant challenges in modern healthcare, making any progress in our understanding of amyloid fibrillization crucial, holding implications for a wide range of debilitating medical conditions.

\begin{figure}[h!]
    \centering
    \includegraphics[width = .85\textwidth]{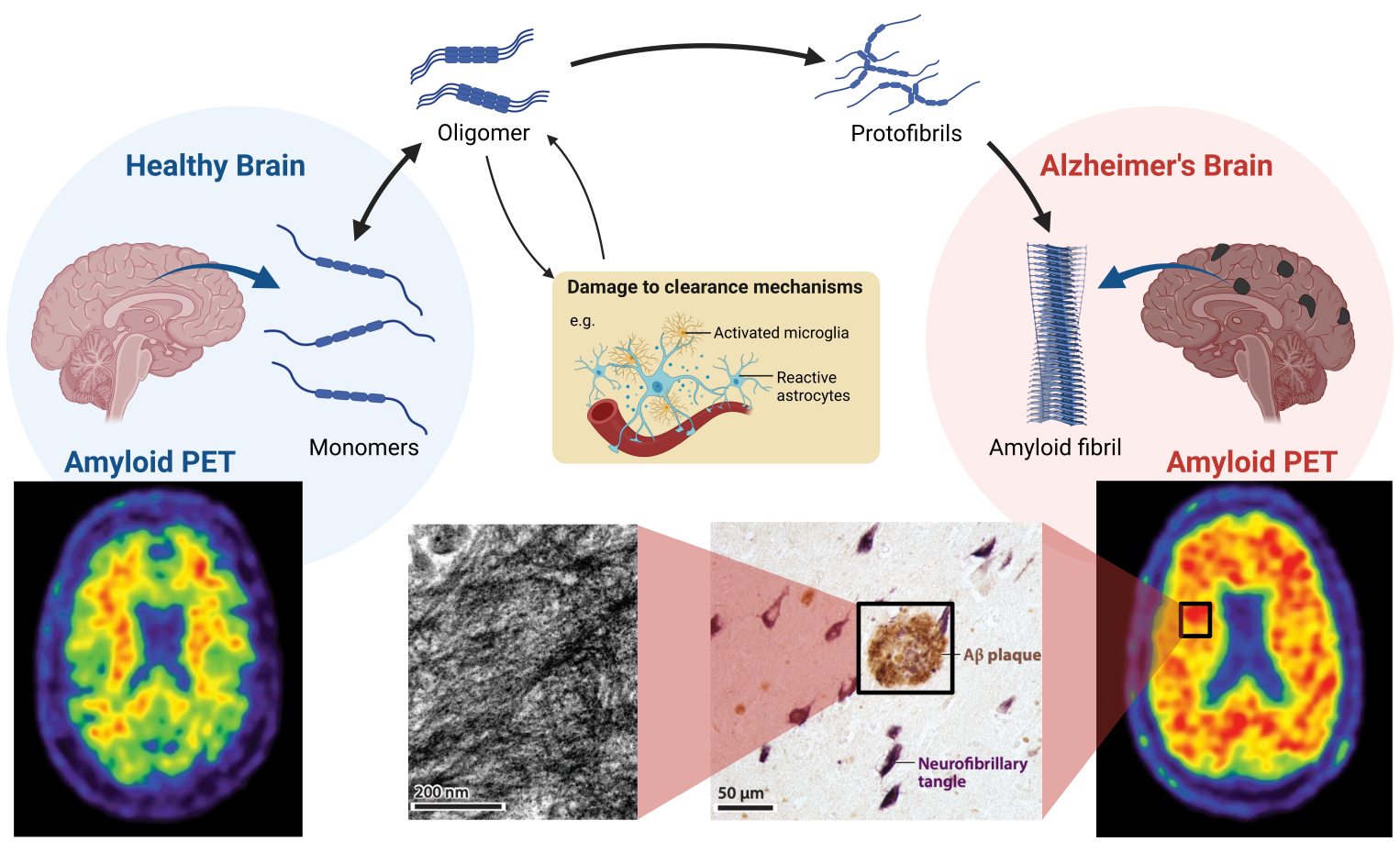}
    \caption{The process of aggregation, from a healthy brain to an unhealthy brain with damaged clearance and saturated in amyloid fibrils, as shown also in amyloid PET images adapted with permission from \cite{ten2018secondary}, and also shown in a light micrograph of A$\beta$ plaques, and an electron micrograph showing a mass of extracellular amyloid fibrils in an A$\beta$ plaque at the microscopic scale, adapted with permission from \cite{walker2015neurodegenerative}.}
    \label{fig:intro1}
\end{figure}

The AD brain has defining pathological features of cortical amyloid plaques, comprising a fibrillar form of the amyloid-$\beta$ protein (A$\beta$) as depicted in \Cref{fig:intro1}. It is thought that such A$\beta$ accumulation facilitates the subsequent cascade of neurofibrillary tangles (NFTs) from aggregated tau protein \cite{lansbury1996reductionist, goedert2017like}. Several hypotheses govern the deterministic modelling of protein accumulation and propagation in AD to gain meaningful estimates of its dynamics. First, the \textit{prion-like hypothesis} \cite{frost2010prion, jucker2018propagation, olsson2018prion, goedert2015alzheimer, mudher2017evidence} broadly postulates that neurodegenerative diseases result from an accumulation of misfolded forms of these proteins, which aggregate and contribute to 
neurodegenerative pathology. In this process, disease-specific misfolded proteins act as a template upon which healthy proteins misfold in a manner akin to prion formation \cite{prusiner1998prions}, forming extensive chains transported through the brain along axonal pathways. Given that aggregates of differing sizes exhibit unique transport characteristics and varying toxicity levels, it is essential to monitor their spatial and temporal evolution independently. 
Second, the \textit{amyloid-$\beta$ hypothesis} \cite{hardy1992alzheimer, hardy1991amyloid,selkoe2016amyloid} establishes that amyloidogenic protein accumulation in AD could be causative, placing A$\beta$ central to disease pathology with supporting experimental evidence suggesting A$\beta$ as a primary driver of AD \cite{selkoe2001alzheimer}. This hypothesis has guided most AD research over the past two decades and motivated the development of many therapeutic antibodies targeting different species of the A$\beta$ peptide \cite{hard2012inhibition}. In the last decade, 15 potential therapeutics targeting the role of A$\beta$ in AD in various ways, including inhibition of enzymes involved in A$\beta$ production and removal of A$\beta$ aggregates using antibodies, have been tested in phase III clinical trials \cite{karran2022amyloid}. However, the failure of most such drug trials 
and recent experimental evidence has renewed scrutiny of its foundational assumptions, arguing for the possible importance of other mechanisms. 
Moreover, an alternative theory, the \textit{A$\beta$ oligomer hypothesis}, suggests that oligomers composed of small numbers of A$\beta$ peptides are the most relevant pathological A$\beta$ species, with amyloid plaques perhaps serving as a reservoir for such species \cite{hong2018diffusible, walsh2020amyloid}. Consistent with this, next-generation therapeutic intervention strategies targeting low molecular weight oligomers of A$\beta$ are showing 
promise \cite{linse2020kinetic}.

Another potential driving factor of AD is the clearance of misfolded proteins, a somewhat elusive and powerful \textit{in vivo} effect which experiments \textit{in vitro} fail to replicate. The production of tau and A$\beta$ peptides is a natural process related to neuronal activity. In a healthy brain, these standard metabolic waste by-products \cite{rumble1989amyloid,bacyinski2017paravascular} are 
removed from intracellular and extracellular compartments by several clearance mechanisms \cite{tarasoff2015clearance, xin2018clearance}. Waste proteins are broken down by enzymes, removed by cellular uptake, crossing the blood–brain barrier, or effluxing to cerebrospinal fluid compartments, eventually reaching arachnoid granulations or the lymphatic vessels. Such healthy clearance mechanisms, working in harmony, avert the buildup of toxic A$\beta$ plaques and tau NFT, but their impairment or dysfunction can lead to AD pathology. The specific descriptions of \textit{in vivo} clearance mechanisms remain a topic of clinical debate; however, the kinetics enabling toxic 
proteins to amass into pathological aggregates can be systematically studied \textit{in vitro} and coupled to dynamic clearance mathematically to simulate the naturally therapeutic, or destabilising, effect of clearance and its response to toxic aggregate mass.

Studies of the history of medicine reveal that significant progress in preventing and treating a disease typically demands a deep understanding of its underlying causes \cite{dobson2013story}. Crucially, however, the molecular mechanisms underlying aggregate proliferation in the complex domain of the brain are still, like clearance mechanisms, poorly understood. The community acknowledges the need for a deeper understanding of molecular processes \textit{in vivo} to achieve success in therapeutic strategies \cite{karran2022amyloid}. Moreover, as emphasised by Karran \textit{et al.} (2022) \cite{karran2022amyloid}, clear experiments of therapeutic hypotheses have been challenging to conduct with anti-A$\beta$ approaches: the target is not clearly defined (amyloid plaques, A$\beta$ oligomers, or monomers), the mechanism by which A$\beta$ affects cognition is unknown (direct or indirect synaptic toxicity, induction of tau pathology, neuroinflammation, or a combination of all and other effects), sometimes failures are attributed to drug administration at too late of a stage of neurodegeneration, and there is also evidence that amyloid removal is faster in patients with high baseline levels further confounding a direct comparison of antibodies \cite{klein2019gantenerumab}. Mathematical modelling can contribute to a deeper understanding of these problems and suggest new approaches. Significant progress have been made in elucidating the molecular mechanisms that occur during the assembly of purified protein molecules under controlled \textit{in vitro} conditions. This advance results from new experimental methods as well as better theoretical models used to analyse the resulting data. Modelling these molecular mechanisms in a disease-relevant system, including \textit{in vivo} effects such as clearance and transport at the brain scale, would provide invaluable insights to guide the design of potential cures for these devastating disorders.

The goal of this study is twofold. The primary focus is to develop and analyse, both analytically and computationally, new models to study protein aggregation kinetics, including \textit{in vivo} effects such as clearance in a brain region locally, then scale up to include transport on the human connectome. Second, we use our new \textit{in vivo} aggregation models to test the impact of potentially therapeutic monoclonal antibodies and inform optimal treatment strategies to maximise toxic mass clearance. Our general approach is to study the size distribution of A$\beta_{42}$ aggregates, with parameters informed by experiments in a HEPES buffer \cite{linse2020kinetic}, both locally and brain-wide using discrete aggregation equations \cite{fornari2019prion, fornari2020spatially, brennanpreprint}, finding analytical relationships where possible and solving numerically on the brain’s connectome. The novelty of our work is exploring how the known \textit{in vitro} aggregate size dynamics and the pharmodynamic effect of drugs \cite{linse2020kinetic,mazer2023development} scales up to affect the progression of neurodegenerative diseases in a more realistic \textit{in vivo} model at the brain scale, exploring a variety of clearance profiles reflective of treatment strategies.
A primary result of interest to the computational biology and pharmaceutical community will be to demonstrate that the non-trivial \textit{in vivo} effects of transport and clearance of oligomers are realized with relatively simple deterministic models and couplings, leading to effects with clear physiological interpretations in neurological disease modelling. Moreover, the mathematical analysis will highlight that clearance mechanisms are crucial in destabilizing the system towards proteopathy and potentially in restabilizing the system, 
with implications for therapeutic intervention within the complex domain of the human AD brain. 

\section{From microscopic models to brain-scale dynamics}
Protein aggregation pathways are complex and involve multiple steps with distinct rates \cite{meisl2017scaling}. Advancements based on chemical kinetics have 
 produced a deep understanding of the fundamental mechanisms underlying the formation of aggregates in ideal conditions at the microscopic scale, thereby making clearer the potential for therapeutic intervention \cite{frankel2019autocatalytic, kundel2018measurement}. Namely, a theoretical framework 
 of the classic nucleated polymerisation theory supplemented with secondary aggregation pathways \cite{cohen2011nucleated1, cohen2011nucleated2, cohen2011nucleated3} has been combined with systematic \textit{in vitro} experiments performed under differing conditions, such as varied concentration or pH \cite{meisl2016quantitative,yang2018role}, to reveal such mechanisms. 
 
\begin{figure}[h!]
    \centering
    \includegraphics[width=\textwidth]{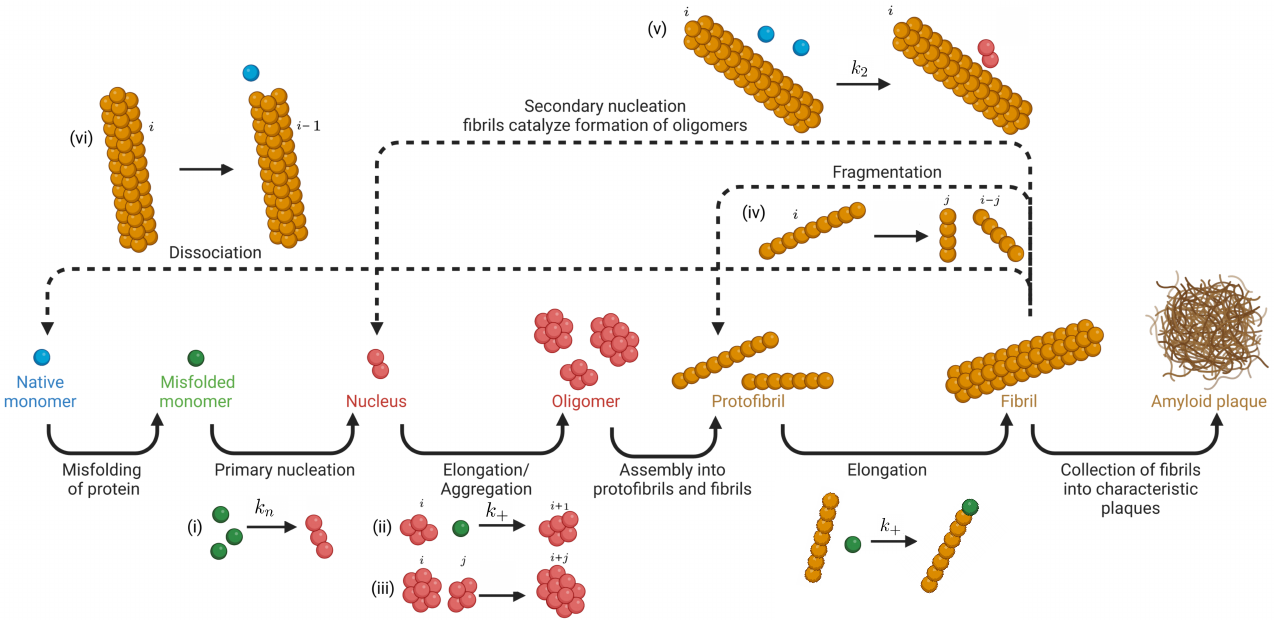}
    \caption{The microscopic nucleation-aggregation-fragmentation processes and secondary pathways. The respective rate constants of processes considered in our models are labelled, as detailed in \Cref{tab:parameters}.}
    \label{fig:intro2}
\end{figure}

 As shown in \Cref{fig:intro2}, it is understood that A$\beta$ monomers misfold by seeding events or prion-like templating \cite{nilsson2002low,jucker2011pathogenic,jucker2013self}, followed by primary nucleation, which leads to the creation of a polymer of length $n_1$ soluble monomers. Filaments elongate and dissociate linearly from both ends with the addition and removal of monomers in a reversible manner. Oligomers and fibrils of different sizes also aggregate to form larger aggregates. Additionally, secondary pathways of fragmentation and monomer-dependent secondary nucleation lead to the creation of new fibril ends from pre-existing polymers.

At present, no brain-scale \textit{in vivo} therapeutic intervention model accounts for the kinetics of protein aggregation. Leading therapeutic models are \textit{in vitro} models for the impact of drugs on the aggregation process \cite{linse2020kinetic}, or brain-wide compartmental models that broadly depict disease development in the presence of drugs over months but neglect the underlying kinetics of drug action \cite{mazer2023development}. 
We combine the aforementioned spatial and time scales, building upon insights from the Smoluchowski-type models of \cite{meislclearance2020} and \cite{fornari2020spatially}, which respectively capture the physiological effects of clearance and brain-wide neuronal transport of aggregates.

Our general approach is to start with a microscopic model with parameters validated experimentally by \textit{in vitro} experiments and extend to both \textit{in vivo} settings and over the entire brain.

First, the  microscopic model that we start with has been carefully calibrated in a series of experiments. In particular, it has been shown that the main processes sufficient to explain the mass increase of oligomers are linear aggregation, primary, and secondary nucleation. Similarly, it has been shown that both depolymerisation and fragmentation can be ignored in the first instance.

Second, we extend the microscopic model to include \textit{in vivo} mechanisms: the natural production of monomers, the clearance of monomers, and a general slowing down due to local cellular effects. This second model represents the local dynamics of oligomers in a given brain region. 

Third, we extend the model to include spatial effects at the brain scale. This is performed by coupling different brain regions and assuming that proteins are mostly transported through the brain along axonal pathways. This final system belongs to the general class of network diffusion-aggregation-fragmentation models discussed in \cite{fornari2020spatially}. In this way, we derive an entirely mechanistic brain-wide network model of neurodegeneration in which we simulate and optimize therapeutic intervention.

\subsection{A minimal microscopic model}
Our minimal microscopic model for protein aggregation includes the following mechanisms: primary heterogeneous
nucleation; secondary nucleation; linear elongation.
Each aggregate of a given size is represented by a population. Let $p_i=p_{i}(t)$ be the concentration of aggregates of size $i = 1,2,3,\dots$. Then, the \textit{microscopic} \textit{master} equations are:
\begin{subequations}
    \label{eq:vitro1}
 \begin{eqnarray}\label{mbeta1}
&&{\frac{{\text{d}}{{m}}}{{\text{d}}t}}  
=-2 k_{n}-2 k_{+}m P-2 k_{2} \sigma(m) m^2 M,\\ 
&&  {\frac{{\text{d}}{{p_2}}}{{\text{d}}t}}  
=k_{n} -2 k_{+}m p_{2}+k_{2} \sigma(m) m^2  M,\\
&&  {\frac{{\text{d}}{{p_i}}}{{\text{d}}t}}  =  %
2 k_{+}m (p_{i-1}-p_{i}),\quad i>2
\label{mbeta5}
\end{eqnarray}
\end{subequations}
where
\begin{equation}
\sigma(m) = \frac{K_m }{K_m + m^{2}}, \quad %
P=\sum_{i=2}^{\infty} p_{i}, \quad  M=\sum_{i=2}^{\infty} i p_{i},
\end{equation}
where $m=p_1$ is the concentration of monomers and the kinetic rate constants are defined in \Cref{tab:parameters}.
Here, $P$ and $M$ are the first two moments 
of the population distribution; they represent the total number and total mass of aggregates, respectively.

\begin{table}[h]
\caption[Kinetic parameters]{Typical  parameters for the  A$\beta$ model. HEPES refers to the buffer 
used for the experiments.}
      \centering
        \begin{tabular}{l|l|c|l}
\hline
param.& mechanism&  A$\beta$42 HEPES \cite{linse2020kinetic}&units\\
 \hline\hline
$k_n$ & primary nucleation      		& $1.6 \times 10^{-11}$ &M\,h$^{-1}$\\
$k_2$ & secondary nucleation & $2.1 \times 10^{14}$&M$^{-2}$h$^{-1}$\\
$K_m$ &monomer saturation      		& $2.3 \times 10^{-17}$&M$^{2}$\\
$K_M$ &mass saturation      		& $2.3 \times 10^{-17}$&M$^{2}$\\

$k_+$ & elongation  			 	&  $1 \times 10^{10}$  &M$^{-1}$h$^{-1}$\\
$m_0$ &initial monomer c. 	& $3 \times 10^{-6}$ &M  \\
 \hline\hline
\end{tabular}
        \label{tab:parameters}
\end{table}

\begin{figure}[ht]
\centering \includegraphics[width=.7\linewidth,trim={1.2cm 0 0 0},clip]{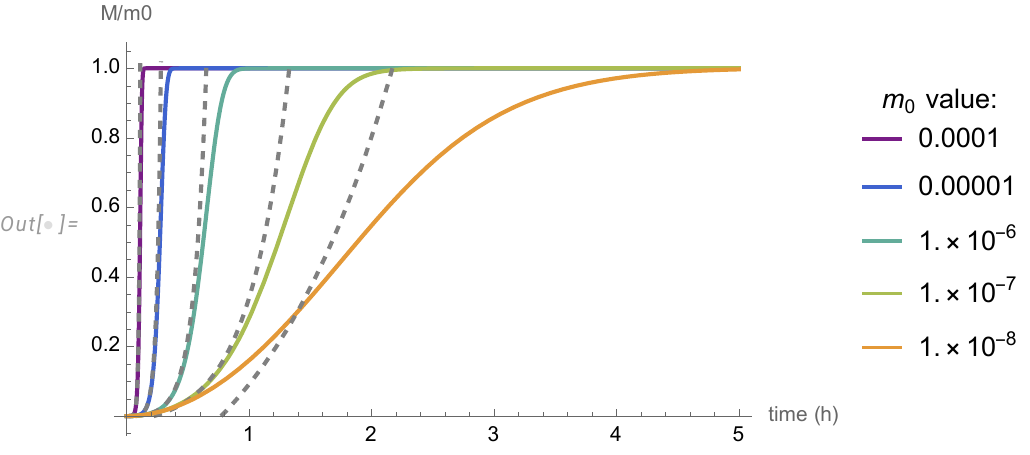}
\caption{Dynamics of the microscopic minimal model based on parameters in \Cref{tab:parameters} for various initial monomer concentrations $m_0$. For large variations of the initial monomer concentration, the dynamics is saturated within hours. The dashed curves show the linear approximation.}
\label{Figkin2}
\end{figure}

It is of interest, before we consider other effects, to understand the overall behavior of this system. This can be easily accomplished by looking at the moment equation obtained as a closed system for $m, M$ and $P$:
\begin{eqnarray}
&&{\frac{{\text{d}}{{P}}}{{\text{d}}t}}  
\,\,=  k_{n}+ k_{2}\, \sigma(m) M,
 \label{PP1}\\
&&{\frac{{\text{d}}{{M}}}{{\text{d}}t}}  
=2 k_{n}+{2 k_{+}m P} + {2}k_{2}\, \sigma(m) M,
 \label{MM1}\\
&&  {\frac{{\text{d}}{{m}}}{{\text{d}}t}}  
\,=-\,2 k_{n}-{2 k_{+}m P}-{2} k_{2}\, \sigma(m) M.
 \label{m}
\end{eqnarray} 
By construction, the total mass of the system is conserved and assuming $M(0)=0$ and $m(0)=m_0$, we have $m_0=M(t) + m(t)$ for all time. As shown in Fig.~\ref{Figkin2}, the behavior of this system is rather simple. The toxic mass $M$ increase from $0$ to $m_0$ in a sigmoid-like manner. Various approximations of this curve have been proposed and used for model validation and parameter fitting \cite{meisl2014differences,frankel2019autocatalytic,cohen2013proliferation}. We observe that for large variations of the initial concentration, the dynamics saturates within hours. 

To better understand the time scale involved in the process, we defined the halftime $\tau$ to be the the time at which an initial unseeded system reaches half of the final concentration $M(\tau)=m_0/2$. To approximate $\tau$, we assume that $k_n$ is small compared to other rates 
and linearize the system~(\ref{PP1}--~\ref{MM1}) with $m(t)=m_0-M(t)$ around $M(0)=P(0)=0$. Expanding $M=\varepsilon \tilde M+\mathcal{O}(\varepsilon^2),P=\varepsilon \tilde P+\mathcal{O}(\varepsilon^2)$,  where $k_n=\varepsilon \tilde k_n$, we obtain, to first order
\begin{eqnarray}
&&{\frac{{\text{d}}{{\tilde P}}}{{\text{d}}t}}  
\,\,=  k_{n}+a \tilde M,\quad {\frac{{\text{d}}{{\tilde M}}}{{\text{d}}t}}  
=2 k_{n}+(b-a) \tilde P + {2}a\tilde M,
 \label{M1}
\end{eqnarray}
where
\begin{equation}
a=\frac{k_2 m_0^2 K_m}{K_m+m_0^2},\qquad b=2 m_0 k_{+}+\frac{k_2 m_0^2 K_m}{K_m+m_0^2}.
\end{equation}
A very good approximation of the exact solution of this linear system with unseeded initial conditions $\tilde M=\tilde P=0$ is obtained by neglecting the fast decaying exponential:
\begin{equation}
M(t)\approx \frac{k_n}{2a} {\left(\left(1+\sqrt{\frac{a}{b}}\right) e^{t(a+\sqrt{ab})
   t}-2\right)}.
\end{equation}
The linear solution gives a reasonable estimate of the initial dynamics and is simple enough to provide an analytical estimate of the halftime:
\begin{equation}\label{taulin}
\tau_{\text{lin}}=\frac{1}{a+\sqrt{a b} }
\log \left(
\frac{\sqrt{b}(a m_0+2  k_n)}{k_n(\sqrt{a}+{\sqrt{b}})}\right).
\end{equation}
In the range of parameters involved, we can further simplify this expression to 
\begin{equation}\label{taulin2}
\tau_{\text{lin}}\approx \frac{1}{ \sqrt{2 k_+ k_2  K_m m_0}}\log \left(\frac{k_2  K_m m_0}{k_n}+2\right).
\end{equation}
\begin{figure}[h!]
\centering \includegraphics[width=.6\linewidth,trim={1.2cm 0 0 0},clip]{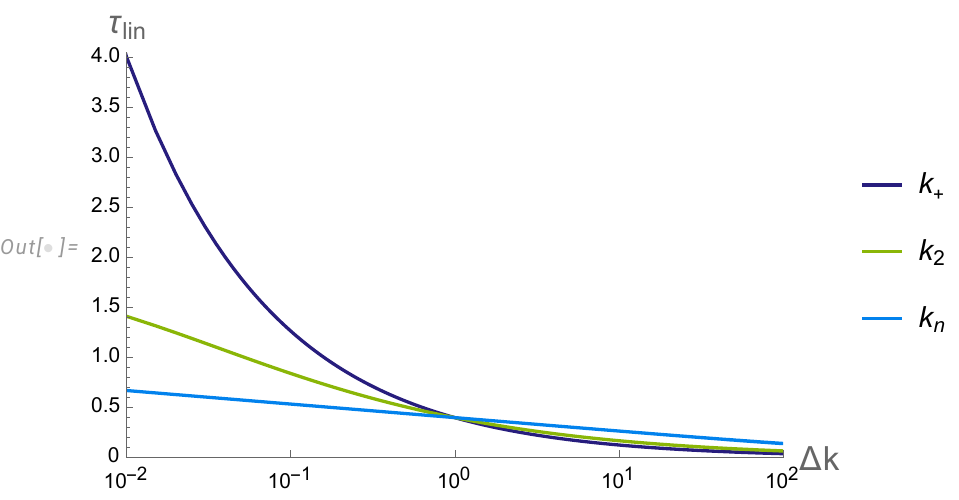}
\caption{Variations of halftimes $\tau_{\text{lin}}$ due to a change of each parameter of the form $k\to (\Delta k) k_{\text{ref}}$, around the reference value $k_{\text{ref}}$ given in \Cref{tab:parameters} (and $m_0=10^{-6}$ M). Note that in this range of parameters, the expressions~(\ref{taulin}) and~(\ref{taulin2}) for $\tau_{\text{lin}}$ are indistinguishable.}
\label{Fig-half}
\end{figure}

We now turn to the variation of time scales associated with the various parameters. We fix the initial concentration to $m_0=10^{-6}$M and systematically vary each parameter over four orders of magnitude around their values given in \Cref{tab:parameters}, to study which rate constant has the most effect on the dynamics. From~(\ref{taulin2}), we see that variations in $m_0$, $k_2$ or $K_m$ are equivalent as these parameters only enter as a product. As expected, the variations due to a reduction in the primary nucleation terms are weak (and only due to the logarithm), whereas reductions in $k_2$ (with a log dependence divided by a root) are faster and variations in $k_+$ dominate (as they depend only on the square root). 

The overall conclusion is that the pure conversion process of a population of monomers into oligomers is governed \textit{in vitro} by rate constants that provide a typical time scale of hours, even when kinetic rates are modified by orders of magnitude. However, we know that \textit{in vivo} any neurodegenerative disease evolves on time scales of years or decades. Even in the case where all kinetic rates are reduced by a factor $T$, the resulting new halftime is simply $\tau T$, which would require a factor of about $10^4$ to reach a time scale of a year. This apparent contradiction requires new mechanisms to explain the dramatic slowdown of these processes.

\subsection{A local model including clearance, production, and saturation}
Based on our understanding of the minimal \textit{in vitro} model that represents a closed-system, we can now extend it to include effects that appear \textit{in vivo}. First, we  assume that there is a regulation mechanism for the production of monomers such that their concentration remains mostly constant, despite their uptake in the formation of larger structures. Therefore, within the modelling framework, we take $m_0$ to be either constant in time or a given function depending on aging and other environmental factor with slow time variation. Second, we  assume that   the main autocatalytic mechanism of secondary nucleation is saturated with respect to the total mass $M$ (rather than $m$ in the initial model). Third, it is believed that clearance is important in both the initiation and evolution  of the disease. Therefore, we assume a linear clearance model with clearance parameters $\lambda_i$ proportional to the oligomer concentration. Assuming, in the first instance, that clearance does not evolve in time, we have
\begin{subequations}
    \label{eq:vivo4}
 \begin{eqnarray}
&&{\frac{{\text{d}}{{m}}}{{\text{d}}t}}  
=0,\\ 
&&  {\frac{{\text{d}}{{p_2}}}{{\text{d}}t}}  
=-\lambda_2 p_2+k_{n}m^2 -2 k_{+}m p_{2}+k_{2} \sigma(M) m^2  M, \label{constm2}\\
&&  {\frac{{\text{d}}{{p_i}}}{{\text{d}}t}}  = -\lambda_i p_i+ %
2 k_{+}m (p_{i-1}-p_{i}),\quad i>2
\label{constm_i}
\end{eqnarray}
\end{subequations}
where
\begin{equation}\label{eqn:moments}
\sigma(M) = \frac{K_M }{K_M + M^{2}}, \quad %
P=\sum_{i=2}^{\infty} p_{i}, \quad  M=\sum_{i=2}^{\infty} i p_{i}.
\end{equation}

We note that in experiments, the process is initiated by seeding with a small amount of oligomers. If dimers are used for seeding, then the initial conditions is simply $m(0)=m_0\ ,P(0)=p_{2}(0), M(0)=2 p_{2}(0)$ and $k_n$ can be taken to be identically vanishing in excellent approximation for $k_n$ sufficiently small. Indeed, once the process is seeded (either through nucleation or oligomer addition), the contribution of the production term $k_n$ becomes negligible in the dynamics. Mathematically, this approximation has the advantage to have $M=P=0$ is a fixed point of the system. Therefore, in the rest of this {section}, we take this point of view and set $k_n=0$.
Further, in the particular case where clearance is size-independent ($\lambda_i=\lambda> 0,\ i>1$), the moment equations read:
\begin{subequations}
    \label{eq:vivo4M}
\begin{eqnarray}
&&{\frac{{\text{d}}{{P}}}{{\text{d}}t}}  
\,\,= -\lambda P +k_{n}m_0^2+ k_{2}\, \sigma(M) m_0^2 M,
 \label{P}\\
&&{\frac{{\text{d}}{{M}}}{{\text{d}}t}}  
=-\lambda M+{2 k_{+}m_0 P} + {2}k_{2}\, \sigma(M) m_0^2 M.
 \label{M}
\end{eqnarray}
\end{subequations}

\begin{figure}[ht!]
  \centering
    \includegraphics[width = 0.65\textwidth, trim={1.5cm 0 0 0}, clip]{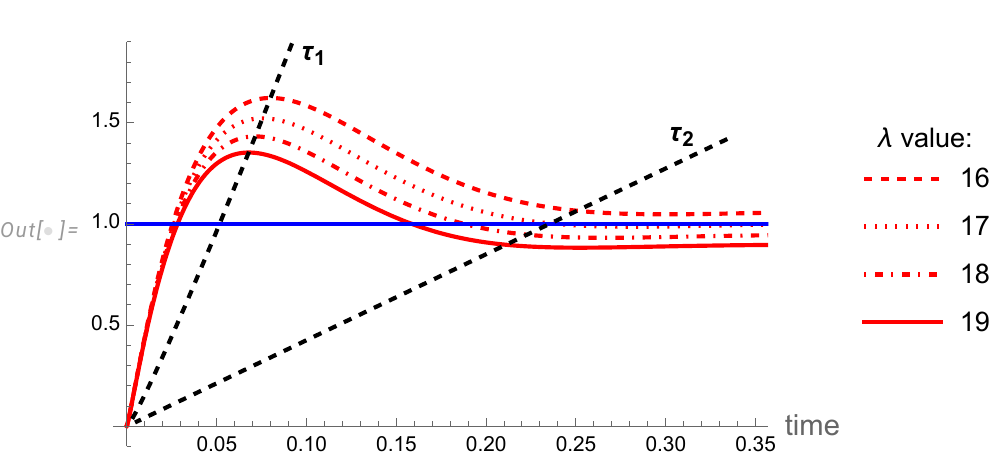}
  \caption{Time evolution of the oligomer toxic ratio $M/m_0$ (red) compared to initial monomer concentration $m_0/m_0$ (blue) in \eqref{eqn:moments}. The system is seeded with $p_2(0) = m_0\times 10^{-4}$. Parameter values are given in \Cref{tab:parameters} with $k_n = 0$. In this case $\lambda_{\text{crit}} = 12705$ from \eqref{eq:agg_crit_cl}. Timescales of interest $\tau_1$ and $\tau_2$ are shown in (b) by the black dashed curves. 
  \label{fig:constant_cl_L0var}}
\end{figure}

This system has two fixed points $P_1=M_1=0$ and
\begin{equation}
\label{eq:fp_agg2}
P_2=\frac{K_M \left(-\lambda ^2+2 \lambda  k_2 m_0^2+2 k_+ k_2 m_0^3\right)}{2 \lambda  \left(k_+
   m_0+\lambda \right)},
   \quad 
   M_2= \frac{K_M \left(-\lambda ^2+2 \lambda  k_2 m_0^2+2 k_+ k_2 m_0^3\right)}{\lambda ^2},
\end{equation}
which exists only if 
\begin{equation}
\label{eq:agg_crit_cl}
0<\lambda<\lambda_{\text{crit}}=k_2 m_0^2+\sqrt{k_2 m_0^3 \left(k_2 m_0+2 k_+\right)}.
\end{equation}
We conclude that this model exhibits a transcritical bifurcation with the property of having the zero trivial state stable for $\lambda>\lambda_{\text{crit}}$ and replaced by a non-vanishing oligomer concentration for $0<\lambda<\lambda_{\text{crit}}$ as shown in \Cref{fig:constant_cl_L0var}. We observe that the critical clearance is independent of the saturation constant $K_M$ and that the asymptotic mass $M_2$ scales with that constant, giving the typical size allowable local oligomer load.

\begin{figure}[h!]
\centering
    \includegraphics[height =.45\linewidth, trim={1.3cm 0 0 0}, clip]{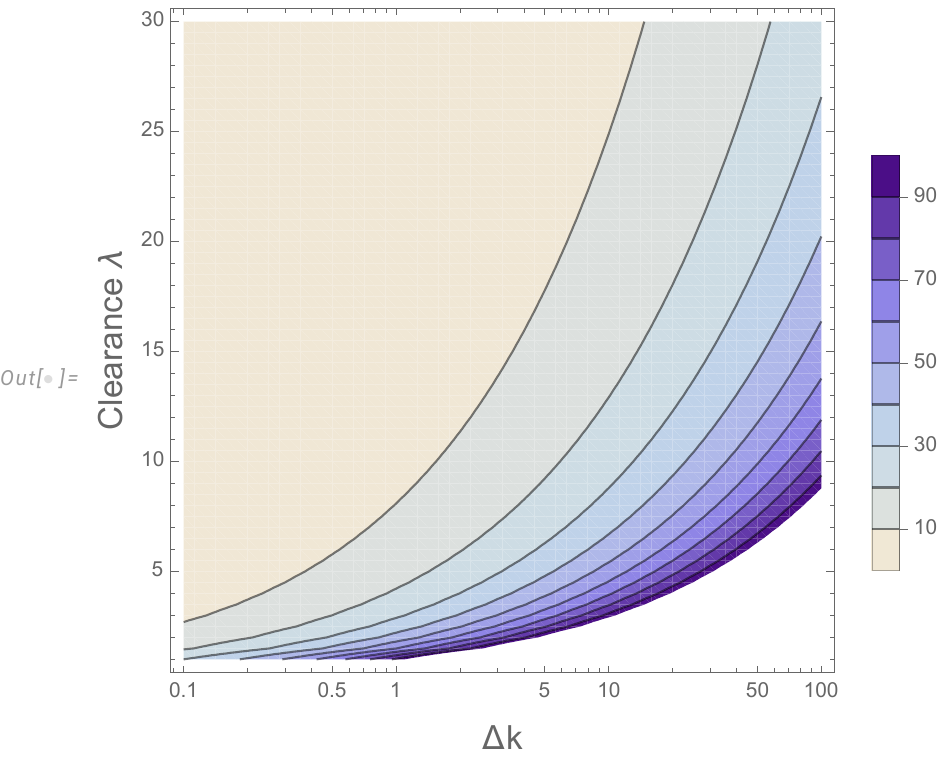}
  \caption{Contour plot of $M_{\text{max}}$ for varying $\lambda$ and $k_2 = k_2 \Delta k $, where $\Delta k$ varies over three orders of magnitude from $10^{-1}$ to $10^2$.  
  \label{fig:contour1}}
\end{figure}

With the introduction of non-zero constant clearance $\lambda$ in the system \eqref{eq:vivo4M}, there are now two timescales associated with the dynamics: {a period of growth up to the maximum toxic mass at time $\tau_1\approx 1/\nu_1$ where $\nu_1 = \sqrt{2k_{+}m_0}-\lambda $ is the first eigenvalue of the system \eqref{eq:vivo4M} when linearized, and a period of decay to time $\tau_2$ when $M-M_2=\epsilon$ where $\epsilon\ll1$} (\Cref{fig:constant_cl_L0var}). Characteristic time scales for the amplification of the aggregate mass are obtained numerically and displayed in \Cref{fig:contour1}. The impact of varying $k_2$ and $\lambda$ on $M_2$ can also be seen by \eqref{eq:fp_agg2}. The dependencies of both timescales on clearance and secondary nucleation are similar. Specifically, we observe that the impact of varying clearance on the timescales and toxic load in the system is profound.

Following the lag period, toxic mass $M$ is dominated by elongation and secondary nucleation, up to time $\tau_1$ when the dampening effects of clearance and saturation of secondary nucleation begin to dominate, initiating a decline in toxic mass from $M_{\text{max}}$ to the steady state $M_2$. The peak toxic mass $M_{\text{max}}$, fixed point $M_2$, and associated timescales $\tau_1$ and $\tau_2$, decrease with increased clearance. Increased secondary nucleation rates $k_2$ has an effect of increasing $M_{\text{max}}$ and $M_2$ and decreasing $\tau_1$ and $\tau_2$, as expected from \eqref{eq:fp_agg2}. Unless $k_2$ is altered dramatically, lowered by one order of magnitude, it has a marginal effect on lowering $\tau_1$ and $\tau_2$ with more effects in lowering $M_{\text{max}}$ and 
$M_2$.

Increasing clearance $\lambda$ has a more significant impact on decreasing toxic load in the local brain system \eqref{eq:vivo4} and both timescales $\tau_1$ and $\tau_2$, and thus dominates in the preferential dampening effects at each phase of the disease cascade. Note that this model does not capture the full disease timescales in a human brain; transport and dynamic clearance have significant effects, as incorporated in the following sections.

If a steady state $p_i^*$ for $i\geq 2$ exists (i.e. if $\lambda_i$ is below some critical clearance threshold and $M$ is bounded), setting \eqref{constm_i} to zero, it must satisfy the recurrence relation
\begin{equation}
\label{fpsoln}
    p_i^* = \delta_i p_{i-1}^*,\quad \delta_i = \frac{2k_+ m_0}{\lambda_i + 2k_+ m_0},\quad i>2,
\end{equation}
where each recursive effect is dependent on $i$. Consequently, the fixed point of the $i$-mer is 
\begin{equation}
\label{analyticalp}
    p_i^* = \triangle_i p_2^*, \quad \text{where } \quad 
    \triangle_i =
\begin{cases}
    \prod_{j=3}^i \delta_j & \text{for } i>2, \\
    1 & \text{for } i = 2.
\end{cases}
\end{equation}
By definition, 
\begin{equation}
\label{Mfp}
    M_* = \sum_{k=1}^\infty k \triangle_k p_2^* = \triangle p_2^*
\end{equation}
where we define $\triangle = \sum_{k=2}^\infty k \triangle_k$. 
Using this in solving \eqref{constm2}, we obtain
\begin{equation}
    p_2^*= \frac{\sqrt{23} \sqrt{-\lambda_2 - 2 k_{+} m_0 + k_2 m_0^2 \triangle}}{1000 \sqrt{
 \lambda_2 \triangle^2 + 2 k_p m_0 \triangle^2}}.
\end{equation}
Thus, a fixed point exists if the following three conditions are satisfied:
\begin{align}
\label{cond1}
    &\text{(C1)}\quad \lim_{{i \to \infty}} \triangle_i = 0,\\
\label{cond2}
    &\text{(C2)}\quad  \triangle = \sum_{k=2}^\infty k \triangle_k \quad \text{converges},\\
\label{cond3}
&\text{(C3)}\quad -\lambda_2 - 2 k_{+} m_0 + k_2 m_0^2 \triangle > 0
\end{align}
where (C1) reflects that the final concentration of aggregates must be decreasing with $i$, (C2) $M_\infty$ must be bounded, and (C3) requires $p_2^*\in \mathbb{R}$. 
Following methods presented in \cite{meislclearance2020} we obtain the following analytical critical clearance formulae dependent on the form of clearance. This case is particularly interesting 
since, in the derivation of steady states, we see that a critical clearance can be established for different forms of size-dependent clearance. Remarkably, our results suggest that, depending on the specific size dependence, 
the processes of elongation and secondary nucleation contribute to the value of the critical clearance to different degrees. An important implication is that, depending on the specific mechanism of clearance, inhibition of aggregation should target different 
processes to reduce the critical clearance rate.

First considering constant clearance, as in \eqref{eq:agg_crit_cl}, (C1--C3) are satisfied for
\begin{equation}
\label{eq:Lcrit}
    \lambda<\lambda_{\text{crit}}=k_2 m_0^2+\sqrt{2 k_2 m_0^3 \left(k_2 m_0+ k_+\right)}.
\end{equation}

\begin{figure}[h!]
\centering \includegraphics[height=.16\textheight,trim={0 0 0 0},clip]{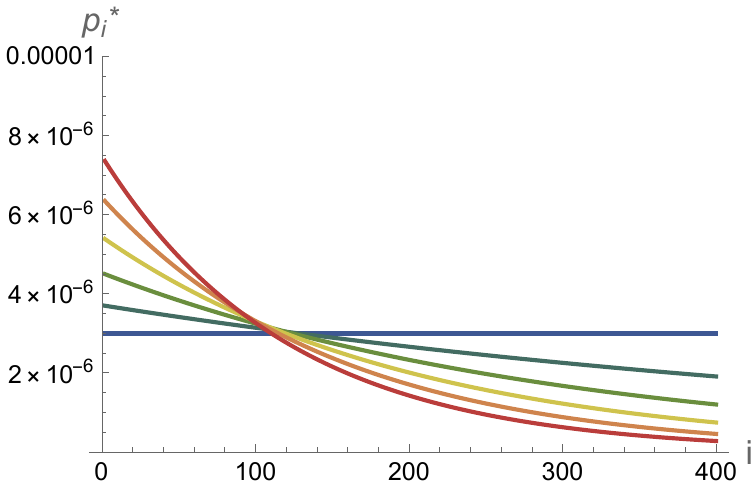}
\includegraphics[height=.16\textheight,trim={0 0 0 0},clip]{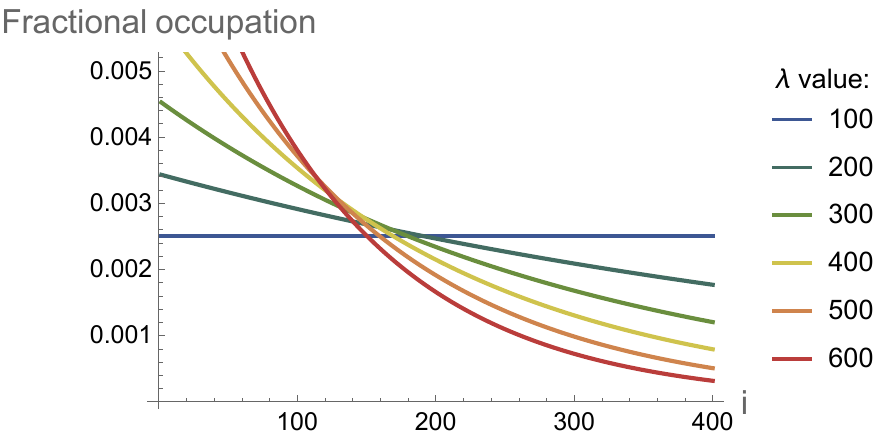}
\caption{The analytical result for the equilibrium length distribution of $p_i^*$ \eqref{analyticalp} in \eqref{eq:vivo4} is shown for various constant clearance rates, with parameters applied from \Cref{tab:parameters}. Heightened clearance has the least impact on oligomers. 
\label{fig:invariantoligomers}}
\end{figure}

The size distribution for aggregates at equilibrium $p_i^*$ for varying subcritical constant clearance rates reveals that the relative occupation of oligomers in the region appears invariant to clearance. Increasing clearance uniformly in the constant clearance regime targets larger aggregates as seen in \Cref{fig:invariantoligomers}. 
This is consistent with Oosawa theory \cite{oosawa1975thermodynamics, oosawa1970size}, which predicts that the length distribution initially develops into a Poisson distribution in the time taken for the monomer-polymer equilibrium to be established before relaxing over a longer time scale 
into an exponential distribution \cite[Figure 5]{cohen2011nucleated3}. Due to the constant supply of monomers, oligomer concentrations are relatively higher, reflective of a brain region.

Considering size-dependent clearance $\lambda_i$, the infinite system \eqref{eq:vivo4} does not yield a closed system for the moments. Still, the fixed point solution \eqref{fpsoln} enables us to determine steady states for aggregate concentrations, subject to conditions for the existence of such a steady state (C1)-(C3). These existence conditions also facilitate the identification of the bifurcation point in clearance, $\lambda_{\text{crit}}$, above which aggregation is negligible. 

First, consider clearance that increases with aggregate size, $\lambda_i = i\lambda_0$. The biological understanding of this case might correspond to clearance mediated by drugs with preferential binding towards larger aggregates. 
In this case,
\begin{equation}
     \lambda_{0,\text{crit}} = \frac{k_{+}k_2 m_0^3 - 1}{k_+ m_0 - k_2 m_0^2}.
\end{equation}
Since $\lambda_i$ is increasing with aggregate size like $\lambda_i = \lambda_0 i$, this intuitively suggests that a lower $\lambda_0$ is required to avoid a diseased state.

For comparison, we consider the opposite case where clearance decreases with aggregate size $\lambda_i = \lambda_0/i$ with the biological interpretation that aggregates become more difficult to clear as they increase in size. Then, an analysis of (C1)-(C3) reveal a critical clearance of
\begin{equation}
\label{eq:lcrit_decreasing_cl}
    \lambda_{0,\text{crit}} = \frac{2a(a + 4k_2 m_0^2)}{a+k_2m_0^2},
\end{equation}
which can be further simplified by again noting that $a = 2 k_{+} m_0 \gg k_2 m_0^2$, so we obtain the approximation $\lambda_{0,\text{crit}} \approx 2a$.

\subsection{A local model including aging and damage}

As the mass of toxic proteins increases, it induces local damage that affects the proper function of the vasculature and all clearance mechanisms \cite{bennett2018tau, carrillo2014, canobbio2015role,michalicova2020tau}. To capture these effects, we let the clearance rates evolve in time up to a lower minimal clearance $\mu_i,\ i>1$. The full system of equations now reads
\begin{subequations}
\label{eq:vivo6}
 \begin{eqnarray}
&&{\frac{{\text{d}}{{m}}}{{\text{d}}t}}  
=0,\\
&&  {\frac{{\text{d}}{{p_2}}}{{\text{d}}t}}  
=-\lambda_2 p_2+k_{n}m^2 -2 k_{+}m p_{2}+k_{2} \sigma(M) m^2  M,\\
&&  {\frac{{\text{d}}{{p_i}}}{{\text{d}}t}}  = -\lambda_i p_i+ %
2 k_{+}m (p_{i-1}-p_{i}),\quad i>2,\\
&&  {\frac{{\text{d}}{{\lambda_i}}}{{\text{d}}t}} =\beta_i M (\mu_i-\lambda_i),
\label{mbeta35}
\end{eqnarray}
\end{subequations}
with moment equations given by 
\begin{subequations}
\begin{align}
\frac{\text{d} P}{\text{d}t} &= -\sum_{i=2}^{\infty} \lambda_i  p_i  + k_n m^2   + k_2\sigma(M) m^2 M,\\
    \frac{\text{d} M}{\text{d}t} &= -\sum_{i=2}^{\infty} \lambda_i i p_i  +2k_n m^2 + 2k_{+} m P + 2k_2\sigma(M) m^2 M.
\end{align}
\end{subequations}
All variables and parameters are given in \Cref{tab:parameters}, and we henceforth consider the value and form of $\mu_i$. The dynamics of this system also display a bifurcation and requires sufficient seeding before it bifurcates to a non-trivial solution where $\lambda_i\to \mu_i$ as explained in \cite{brennanpreprint}.
\begin{figure}[htbp!]
  \centering
  \subfloat[$\lambda(0)>\mu>\mu_{\text{crit}}$ dynamics, $\mu = 13000$, unseeded. 
  ]{
    \includegraphics[height =.17\textheight, trim={1.3cm 0 0 0}, clip]{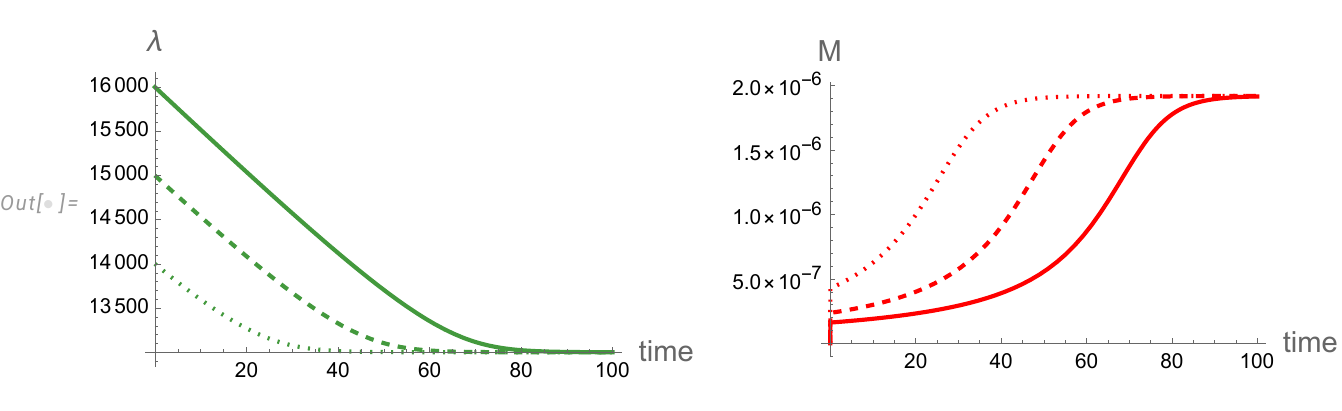}
    \includegraphics[height =.17\textheight, trim={14.5cm 0 0 0}, clip]{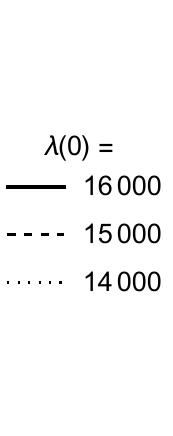}
  }
  \\
  \subfloat[$\lambda(0)>\mu_{\text{crit}}>\mu$ dynamics, $\mu = 9000$, unseeded.]{
    \includegraphics[height =.17\textheight, trim={1.3cm 0 0 0}, clip]{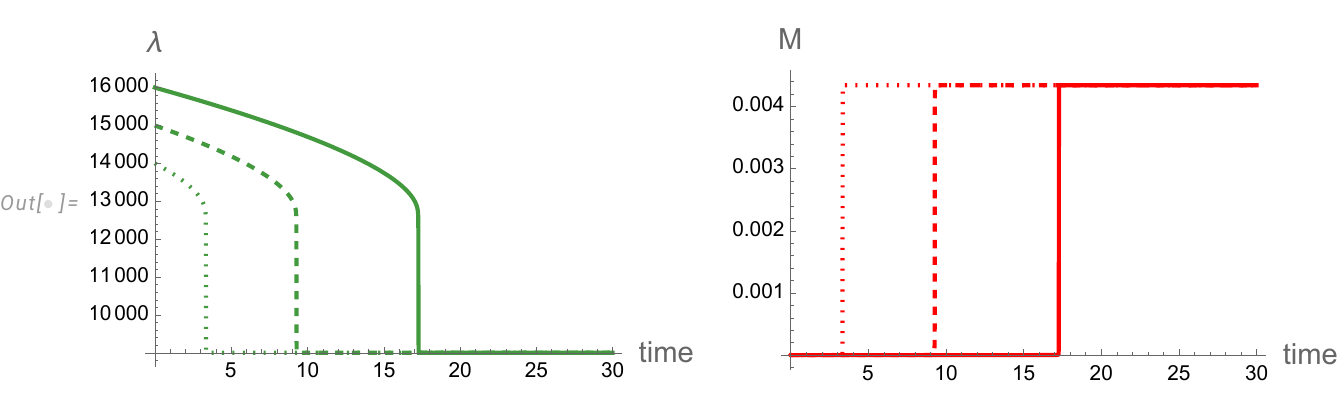}
    \includegraphics[height =.17\textheight, trim={14.5cm 0 0 0}, clip]{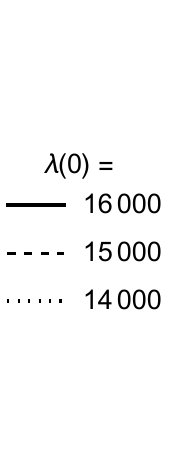}
  }
  \\
    \subfloat[$\mu_{\text{crit}}>\lambda(0)>\mu$ dynamics, $\mu = 9000$, unseeded]{
    \includegraphics[height =.17\textheight, trim={1.3cm 0 0 0}, clip]{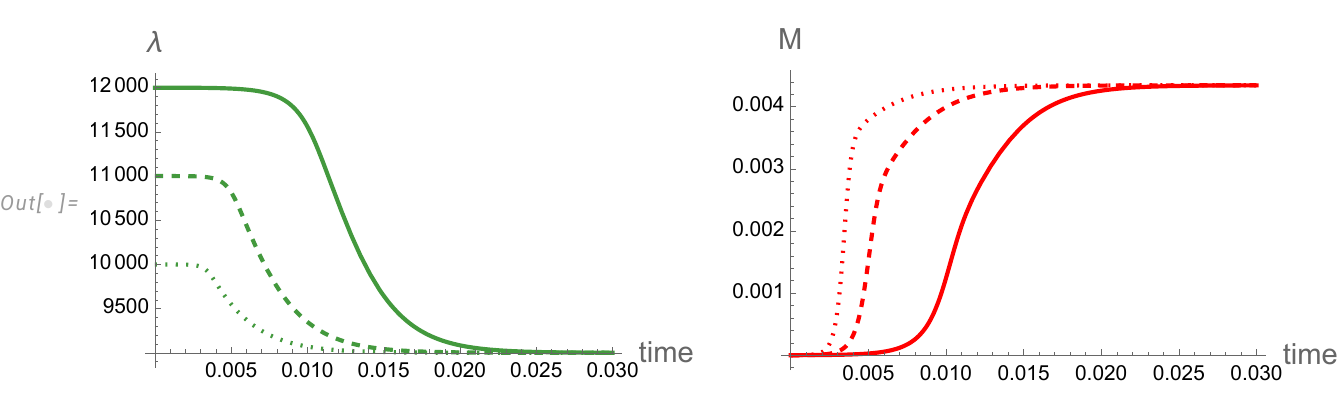}
    \includegraphics[height =.17\textheight, trim={14.5cm 0 0 0}, clip]{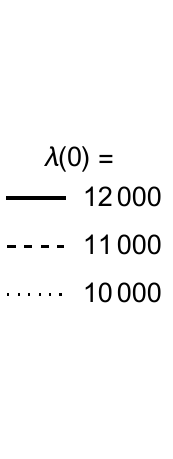}
    }
    \\
    \subfloat[$\mu_{\text{crit}}>\lambda(0)>\mu$ dynamics, $\mu = 9000$, seeded with $M(0)=2\times 10^{-4}$.]{
    \includegraphics[height =.17\textheight, trim={1.3cm 0 0 0}, clip]{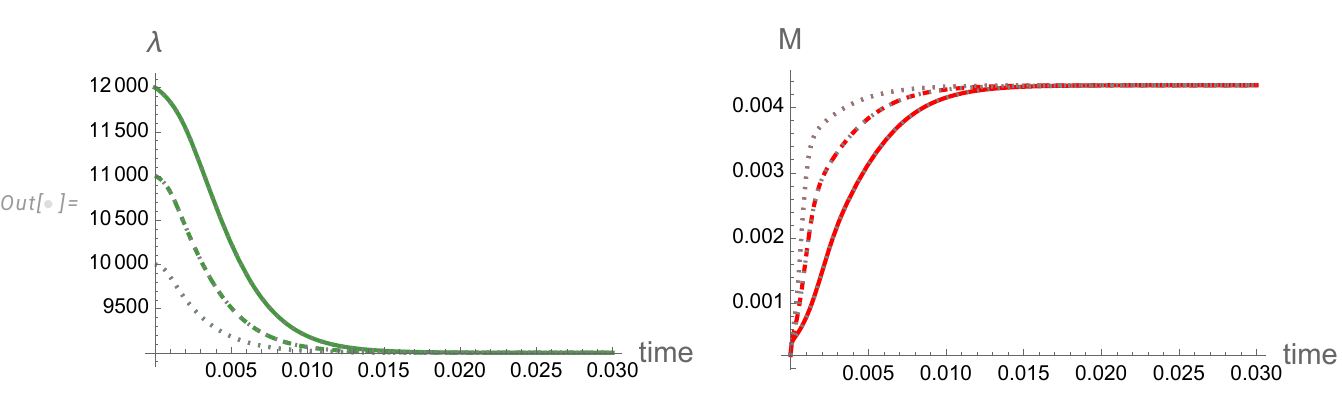}
    \includegraphics[height =.17\textheight, trim={14.5cm 0 0 0}, clip]{figures/legend3.pdf}
  }
  \caption{Toxic mass and clearance dynamics of \eqref{eq:vivo6} corresponding to different uniform initial clearance $\lambda(0)$ and basal clearance capacity $\mu$ values across aggregate sizes. Different parameter regimes are displayed relative to the critical basal clearance capacity $\mu_{\text{crit}}$, with the dynamics of each having distinct timescales of invasion. In all cases, $\mu_{\text{crit}}=12705$ according to \eqref{eq:agg_crit_cl}. Alongside parameter values in \Cref{tab:parameters}, $\beta = 10^{11}$ and concentrations are rescaled by $c = 10^6$. Here, we do not neglect the important effects of primary nucleation, as they are instrumental in slowly decaying clearance. The $k_n=0$ approximation does not approximate the dynamics of (a) or (b) well when seeded, but it does approximate the system well if $\lambda(0)<\mu_{\text{crit}}$, shown by the dotted curves in (d).
  \label{fig:dynamicclconst}}
\end{figure}

If we remove clearance in this model, it becomes equivalent to the previous (constant clearance) system \eqref{eq:vivo4} leading to unbounded growth. Once clearance decays to the basal clearance capacity associated with maximal damage, the dynamics of the system align with those of \eqref{eq:vivo4}. However, the system dynamics prior to the point at which $\lambda_i = \mu_i$ are markedly different, as we observe sigmoidal growth in toxic mass as local clearance is driven to the basal clearance capacity.

\subsection{Coupling the microscopic model to transport}

So far, we have developed a model suitable to describe the dynamics of oligomer concentrations within a single region. Next, we consider a system with multiple regions of interest connected through a network with diffusive transport between different regions. The brain connectome is defined  as a weighted graph ${\mathcal{G}}$ with $V$ nodes ($V$ for  vertices) and $E$ edges obtained from  tractography of diffusion tensor images. From the tractography, we extract the weighted adjacency matrix $\mathbf{A}$ and define the graph Laplacian $\mathbf{L}$
\begin{equation}
L_{ij} = - A_{ij}+\sum_{j=1}^{V} A_{ij},\qquad i,j=1,\ldots,V.
\end{equation}
There are other possible definitions of the graph Laplacian obtained by normalizing rows, columns, or both. However, this is the only graph Laplacian that has the property preserves mass during transport and the requirement that no transport takes place between two regions with the same concentration (for a detailed description and discussion see \cite{brennan_goriely_chapter}).

Incorporating the assumptions specified in the previous section, $k_n=0$ to approximate the seeded system, the truncated master equations for aggregates of size $i$ at node $j=1,2,\cdots,V$, including axonal transport, form an infinite system of first order ODES:
\begin{subequations}
\label{super_general_system}
\begin{align}
&{\frac{{\text{d}}{{m_j}}}{{\text{d}}t}}  
=- \rho_{1}\sum_{k=1}^V L_{jk} {m_{k}},
 \label{mj}\\
&  {\frac{{\text{d}}{{p_{2,j}}}}{{\text{d}}t}}  
=- \rho_{2}\sum_{k=1}^V L_{jk} {p_{2,k}}-\lambda_{2,j} p_{2,j} + k_n m_j^2 
-2 k_{+}m_j p_{2,j}+k_{2} \sigma(M_j) m_j^2  M_j,\\
&  {\frac{{\text{d}}{{p_{i,j}}}}{{\text{d}}t}}  =- \rho_{i}\sum_{k=1}^VL_{jk} {p_{i,k}} -\lambda_{i,j} p_{i,j}+ %
2 k_{+}m_j (p_{i-1,j}-p_{i,j}), \\
&  {\frac{{\text{d}}{{\lambda_{i,j}}}}{{\text{d}}t}} = f(\lambda_{i,j}), \label{eq:super_cl}
\end{align}
\end{subequations}
with initial conditions
\begin{equation}
    m_i(0) = m_{i,0}, \quad p_{i,j}(0)=p_{i,j,0}.
    \label{eq:ics}
\end{equation}
Here $\rho_i$ is the diffusion coefficient of the $i$-mer taken to be small or a function of $i$ \cite{fornari2020spatially,bertsch2017alzheimer}. As before, the first two moments of the population distribution are 
\begin{equation}
    P_j=\sum_{i=2}^{\infty} p_{i,j},\quad M_j=\sum_{i=2}^{\infty} i p_{i,j},
\end{equation}
representing the total number and total mass concentrations of aggregates in the ROI $j$, respectively. 
Taking clearance to be constant in time, the moment equations are given by:
\begin{eqnarray}
\label{eq:brainwide_M}
&&{\frac{{\text{d}}{{M_j}}}{{\text{d}}t}}  
=-\sum_{i=2}^\infty \lambda_{i,j} i p_{i,j}- \rho_{j} \sum_{k=1}^VL_{jk} {M_{k}} 
+2 k_{n} m_j^2 
+2 k_{+}m_j P_j+2 k_{2} \sigma(M_j) m_j^2  M_j, \nonumber \\
 \label{Mj}\\ 
 &&{\frac{{\text{d}}{{P_j}}}{{\text{d}}t}}  
= -\sum_{i=2}^\infty \lambda_{i,j} p_{i,j}- \rho_{j} \sum_{k=1}^VL_{jk} {P_{k}} 
+ k_{n} m_j^2 
+ k_{2} \sigma(M_j) m_j^2  M_j.
\end{eqnarray}

\begin{figure}[h!]
    \centering
    \includegraphics[width = .75\textwidth,trim={0 0 0 0}, clip]{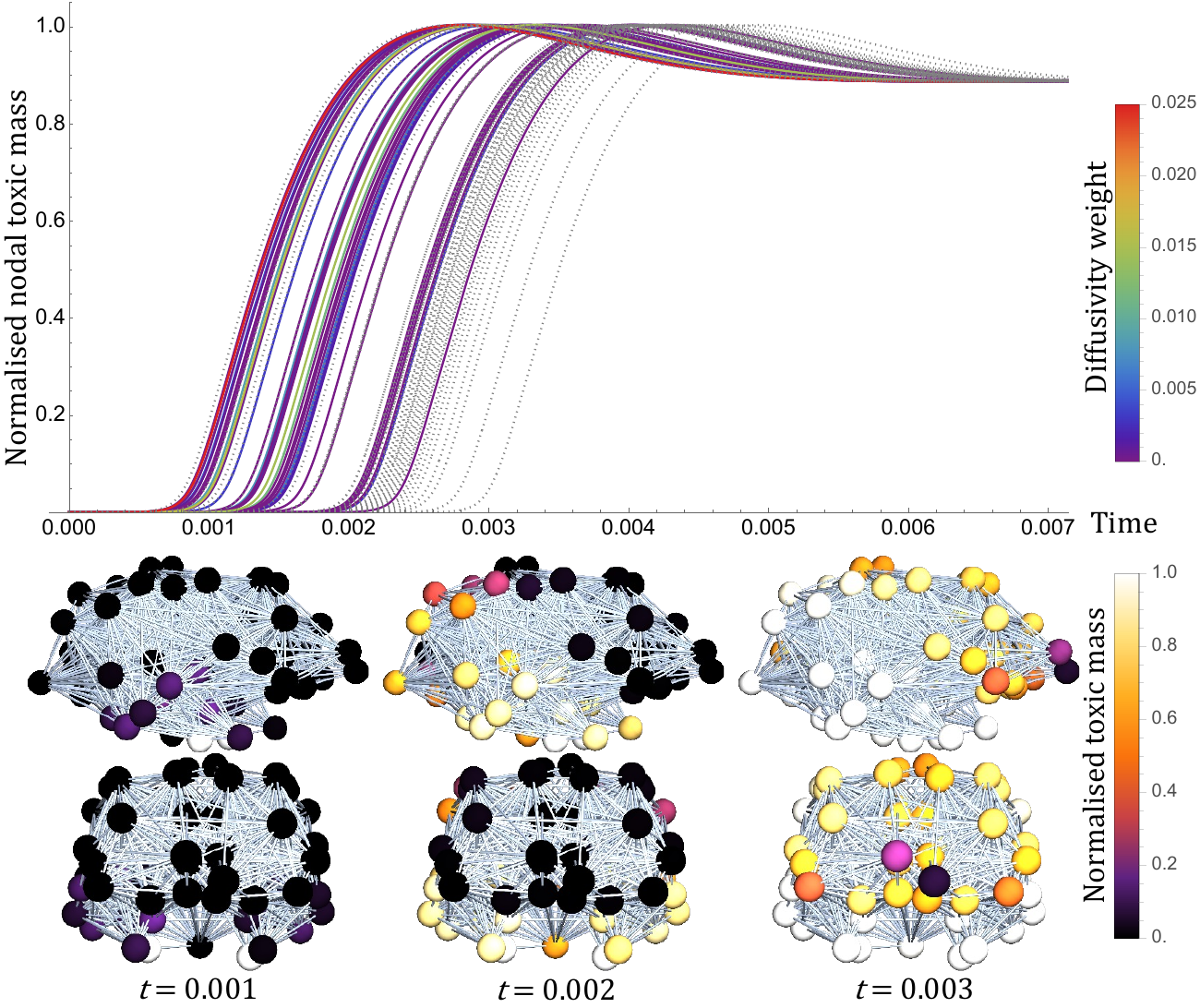}
    \caption{The order of tau invasion of nodes on the connectome and the corresponding sagittal and coronal views of the spread of tau on the brain network in the invasion window $t=0.001$ to $t=0.003$. We seed in the entorhinal cortex $p_2 = m_0$ and display the resulting normalised total toxic mass $M(t)$ on each node, given constant $\lambda = 10^3$ in \eqref{super_general_system}. The top row displays the cascade of toxic proteins through each node, where rainbow colours show the nodal toxic mass in immediate neighbours of the seed node and their diffusivity weighting, and the grey curves correspond to nodes of path length two from the seed node.}
    \label{fig:invasion_order_brainwide}
\end{figure}

The order of node invasion across the brain network suggests a direct applicability of the diagonalization results of \cite{fornari2020spatially} to our model. In \Cref{fig:invasion_order_brainwide}, we observe immediate neighbors of the single seed node are invaded first, influenced heavily by the diffusivity weighting along the edges emanating from this seed. Subsequently, nodes of path length two are invaded, encompassing all nodes in the small-world network. Nodes with the lowest connectivity, such as the frontal pole, are invaded last. 
\Cref{fig:invasion_order_brainwide} also displays the order of ROI invasion across the connectome and the corresponding network representation. 

\begin{figure}[h!]
    \centering
    \includegraphics[width = .9\textwidth]{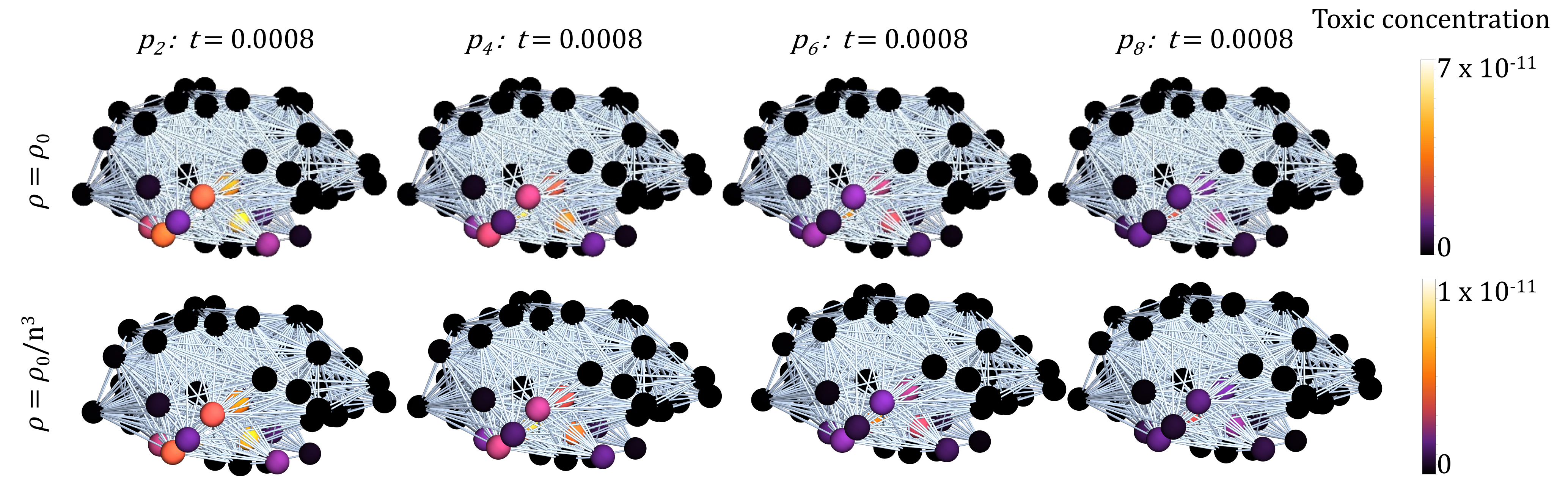}
    \caption{Sagittal views of the distribution of aggregates of increasing size on the brain network at the fixed time point $t=0.0008$. Here we seed in the entorhinal cortex $p_2 = m_0$, with constant clearance $\lambda = 10^3$ in \eqref{super_general_system}. The difference between the top row and bottom row is constant diffusion $\rho=0.01$, and size-dependent diffusion $\rho = \rho_0/n^3$, where $n$ is the oligomer size, resulting in lower toxic concentrations at $t=0.0008$ (see colour bar limits), but no change in the relative occupancy of aggregates.}
    \label{fig:brainwide_figs_progression_aggs}
\end{figure}

The growth of toxic mass across the connectome in this order is facilitated primarily by the transport of nuclei. To demonstrate this, first the distribution of oligomers of size $i=$ 2, 4, 6 and 8, at a fixed time point, is seen in \Cref{fig:brainwide_figs_progression_aggs} (first row) given a size-independent diffusion coefficient in \eqref{super_general_system}. Aggregates spreading in travelling waves, and total toxic mass is a propagating front as described in \cite{putra2023front}. It is observed experimentally \cite{nicholson2000diffusion, goodhill1997diffusion,nicholson1998extracellular} that \textit{in vivo} aggregate transport scales with size. Specifically, large fibril assemblies do not move, and in previous models \cite{bertsch2017alzheimer,achdou2013qualitative} the diffusion coefficient of a soluble peptide is taken to scale as a reciprocal of the cube of its molecular weight. 

Perhaps counter-intuitively, the distinction in progression through nodes being characterised by first low weight oligomers, and later by increasing aggregate sizes, is not exacerbated by the assumption of transport scaling with size, as seen by comparing rows in \Cref{fig:brainwide_figs_progression_aggs}. The order of invasion and relative times of invasion remain the same, but the time taken to reach higher concentrations increases. This emphasises the role of nuclei spreading through the connectome, since local dynamics dominate over the delivery of aggregates by diffusion. We conclude that regardless of the form of size-dependent diffusivity, provided diffusion is small, the transport of larger aggregates does not significantly change the brain-wide dynamics.

\section{Applications to therapeutic modelling}
Next, we investigate the impact of monoclonal antibodies (mAbs), a form of immunotherapy known to influence microscopic parameters in chemical kinetic models, on the aggregation and propagation of misfolded toxic proteins. Using the local and brain scale models derived and studied in the previous section, we analyse, replicate, and propose treatment strategies on the structural connectome based on results \textit{in vitro}. Parameters extracted \textit{in vitro} are scaled up to display the impacts on whole-brain neurodegeneration to simulate treatments in a fully mechanistic model at the whole-brain scale. Notably, FDA-approved drugs such as lecanemab, and drugs with promising phase III trial results, such as donanemab, operate by ultimately increasing the effectiveness of brain clearance. Thus, a model that captures the intricacies of the brain's endogenous clearance and the enhancement of this mechanism from drugs is crucial.

\begin{figure}[h!]
    \centering
    \includegraphics[width = .85\textwidth, trim={1.5cm 0 0 0}, clip]{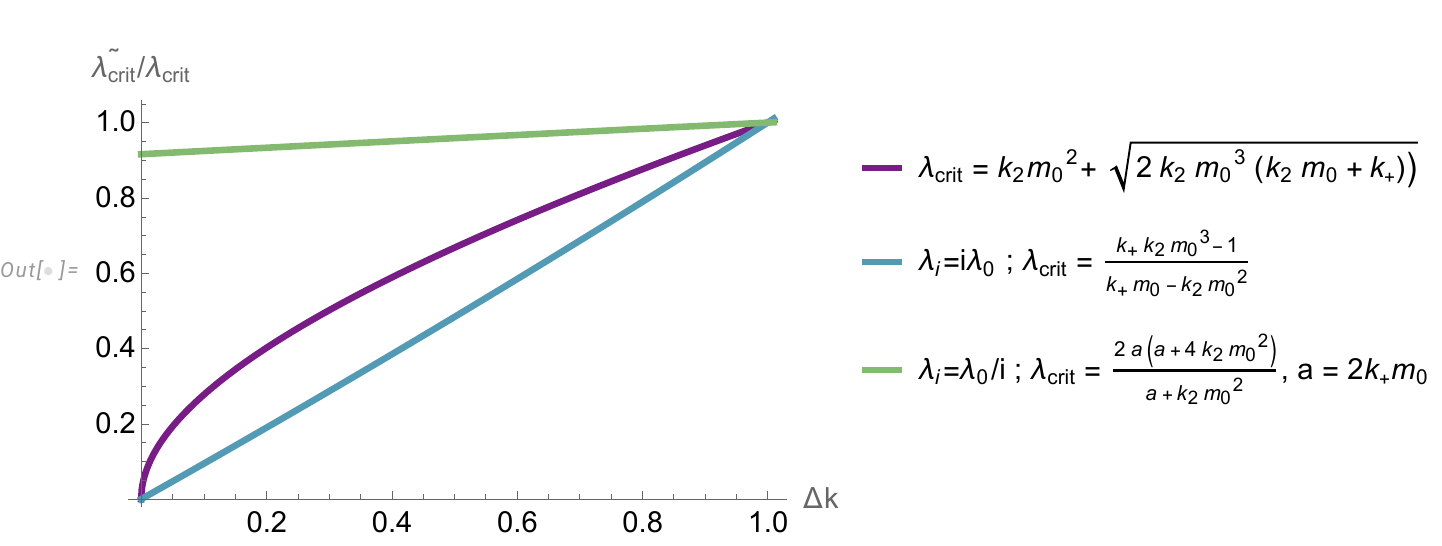}
    \caption{For varying $\tilde{k_2} = \Delta k k_2$, the ratio of critical clearance values in the presence of drugs $\tilde{\lambda_{\text{crit}}}$ and critical clearance associated with $\Delta k = 1$ (without drugs). Results are displayed for critical clearance formulae associated with different models, namely  \eqref{eq:vivo4} with constant clearance $\lambda_i = \lambda$ (purple), and size-dependent clearance $\lambda_i = i \lambda_0$ (blue), and $\lambda_i = \lambda_0/i$ (green). 
    \label{fig:vark2}}
\end{figure}

Linse \textit{et al.} (2020) \cite{linse2020kinetic} found, by fixed point solution method of the coarse-grained protein kinetic equations, the kinetic fingerprints of various mAbs. For example, using the method described in \cite{meisl2016molecular}, they identified that aducanumab causes an apparent reduction of the free oligomer concentration by inhibiting the critical molecular process, secondary nucleation, through which oligomers form. With the highest dose of aducanumab (100 nM), corresponding to a substoichiometric molar ratio of 0.03:1 antibody: A$\beta_{42}$, there is a 69\% reduction in the secondary nucleation rate constant. With the lowest dose tested (250 pM), there is still a 39\% reduction. Thus one effect of the drug on the aggregation pathway can be expressed as a rate constant change in the presence of aducanumab:
\begin{equation}
    \tilde{k_2} = \Delta k k_2,
\end{equation}
where $\Delta k \in[0,1]$ is dependent on the dose of the drug. Importantly, variations in $k_2$ reduce the critical clearance for all models, as seen in \Cref{fig:vark2}, with implications for disease timescales and propagation patterns as discussed in \cite{brennanpreprint}.

The inhibitory effect on secondary nucleation originates from the interaction of aducanumab with amyloid fibrils. Linse \textit{et al.} (2020) \cite{linse2020kinetic} observed a low affinity for monomers and a very high affinity (1nM) for fibrils, in agreement with previous findings \cite{arndt2018structural}. Fibrils become fully coated with aducanumab along their entire length, effectively interfering with secondary nucleation at the fibril surface. Further, mAbs binding to amyloid fibrils targets them for microglia-mediated removal and enhances the clearance of plaques \cite{sevigny2016antibody,soderberg2023lecanemab}. 

Mazer \textit{et al.} (2023) \cite{mazer2023development} model the effect of drugs solely through increased clearance in a compartmental model; the Q-ATN model. The pharmacodynamic (PD) drug effect of antibody-mediated plaque removal is quantified by a linear relationship between antibody concentration and clearance. 
 The PD model assumes a pseudo-first-order rate constant for plaque removal ($\lambda_{\text{drug}}$) that is proportional to the plasma concentration ($C_p$), with an antibody-specific proportionality constant ($L$) given by
\begin{equation}
\label{eq:pharam_clearance}
    \lambda_{\text{drug}}(C_p) = L C_p.
\end{equation}
The $L$ values were estimated by fitting the time course of mean amyloid PET data during treatment using a non-linear least squares algorithm.

We directly substitute drug inhibited parameters into the system \eqref{eq:vivo4} to obtain, in the case of aducanumab,
\begin{subequations}
\label{eq:vivo4_adu}
\begin{align}
& \frac{{\text{d}}{{m}}}{{\text{d}}t} = 0, \\
& \frac{{\text{d}}{{p_2}}}{{\text{d}}t} = -(\lambda_2+\lambda_{\text{drug}}) p_2 + k_{n} m^2 - 2 k_{+} m p_{2} + \tilde{k}_{2} \sigma(M) m^2 M, \\
& \frac{{\text{d}}{{p_i}}}{{\text{d}}t} = -(\lambda_i+\lambda_{\text{drug}}) p_i + 2 k_{+} m (p_{i-1} - p_{i}), \quad i > 2
\end{align}
\end{subequations}
where, as before,
\begin{equation}
\sigma(M) = \frac{K_M }{K_M + M^{2}}, \quad %
P=\sum_{i=2}^{\infty} p_{i}, \quad  M=\sum_{i=2}^{\infty} i p_{i}.
\end{equation}
\begin{figure}[h!]
    \centering
    \includegraphics[width = .45\textwidth]{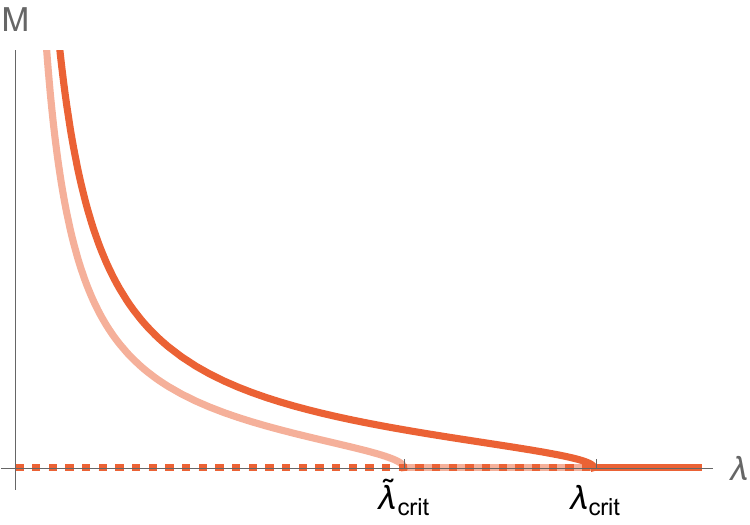}
    \caption{Unstable (dashed) and stable (solid) equilibrium solutions for total toxic mass $M_\infty$ are depicted without drugs (red) and in the presence of aducanumab (lighter red). The presence of aducanumab initiates a reduction in secondary nucleation to $\tilde{k_2} = 0.5 k_2$, causing the critical clearance reduces from $\lambda_{\text{crit}} = 12750$ to $\tilde{\lambda}_{\text{crit}} = 8750$ 
    given by \eqref{eq:agg_crit_cl}, and a reduction in the saturated toxic state $M_2$ given by \eqref{eq:fp_agg2}. 
    }
    \label{fig:toxic_equi}
\end{figure}

Next, we consider the effects of aducanumab on only the kinetic parameters expressed through the reduced secondary nucleation rate constant $\tilde{k_2}$, identified to be the drug's inhibitory action on the aggregation chain. This results in a new critical clearance, which varies with $k_2$ according to \eqref{eq:agg_crit_cl} (\Cref{fig:vark2}), and the steepness of the $M_2$ fixed state according to \eqref{eq:fp_agg2}. The toxic mass equilibrium is shown in \Cref{fig:toxic_equi} in the presence of the kinetic action of aducanumab (orange) and without treatment (red). In addition to this kinetic effect, \cite{mazer2023development} suggest that the preferential binding of aducanumab to fibrils facilitates clearance at a rate proportional to drug concentration. Therefore, we model this effect by increasing aggregate size-specific clearance by $\lambda_{\text{drug}}$. 

\begin{figure}[h!]
    \centering
    \includegraphics[width=\textwidth]{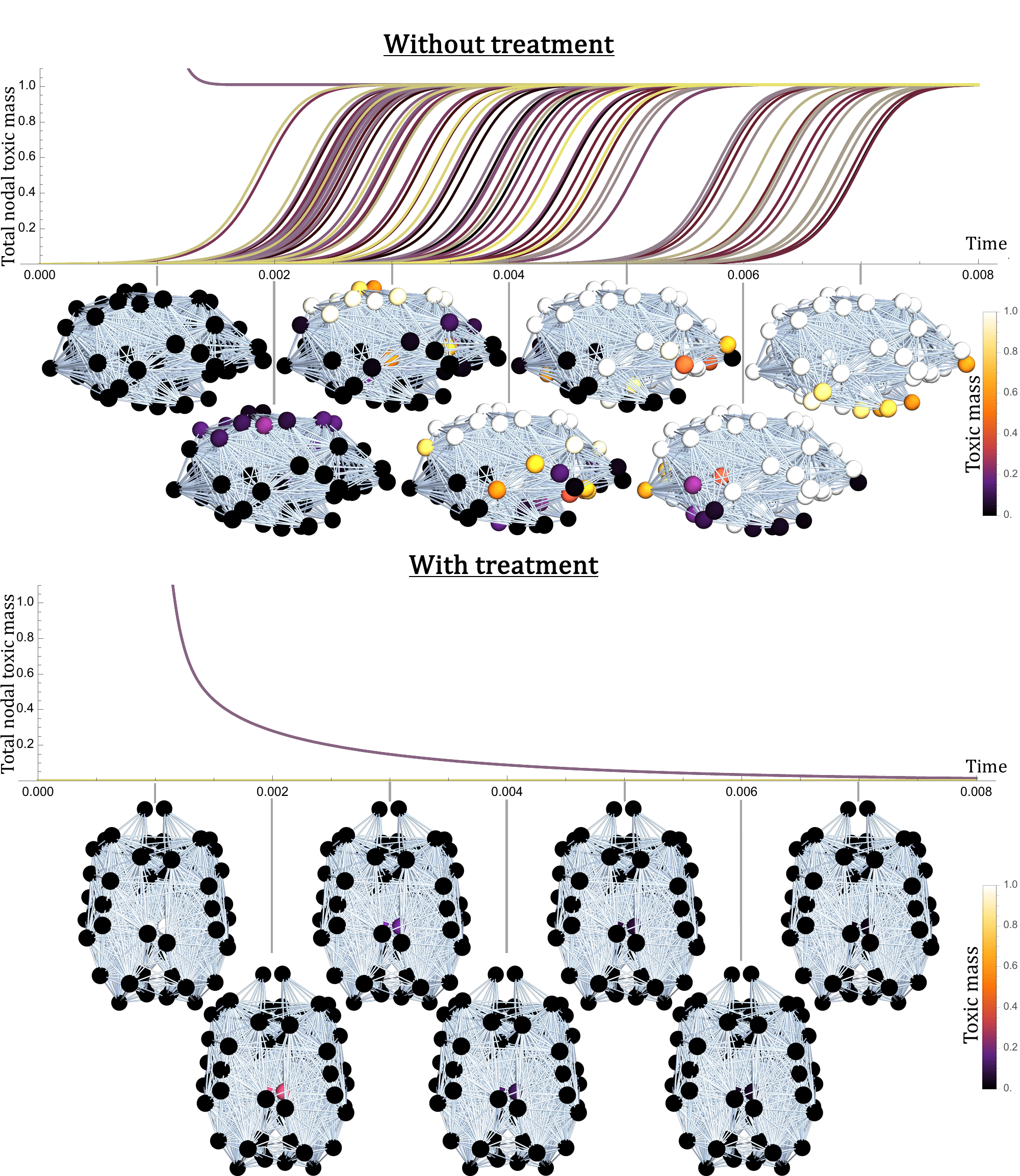}
    \caption{A$\beta$ toxic mass progression through the connectome without treatment (top) and given aducanumab (bottom), initiating a drop in secondary nucleation $\tilde{k_2} = 0.5 k_2$ (without considering the binding effects of enhanced clearance). Subcritical initial clearance of $\lambda(0)=10$ is uniform across nodes, and we impose an initial seed of $M(0) = m_0$ in the posterior cingulate, where A$\beta$ is first observed. Here in \eqref{super_general_system}, 
    $\rho=0.01$, $k_n=0$, and all other parameters are shown in \Cref{tab:parameters}.
    \label{fig:Abetadrugvsnodrug}}
\end{figure}

Next, we couple \eqref{eq:vivo4_adu} to transport across the connectome using the graph Laplacian, as done in \eqref{super_general_system}. In the growth-dominated regime, the homogeneous system is a good approximation of the brain-scale system. The progression throughout the brain network, with and without treatment, is illustrated in \Cref{fig:Abetadrugvsnodrug} for A$\beta$, considering an initial clearance of $\lambda = 9000$ in the network aggregation-clearance model \eqref{super_general_system}. This brain-wide drug simulation highlights the capability of monoclonal antibodies to localise the disease if caught within the clearance window of $\tilde{\lambda}_{\text{crit}}< \lambda<\lambda_{\text{crit}}$. This dependence on drug efficacy on initial clearance, along with the characterisation of the drugs effects purely in terms of brain clearance, further emphasises the pivotal role of clearance in neurodegeneration.

\clearpage
\section{Mathematically informed treatment strategies}

Mathematical models provide a platform for experiments that would be otherwise difficult or impossible to conduct in humans. Here we study potential treatment strategies based on our network model of brain-scale aggregation with \textit{in vivo} effects.

\subsection{Target clearance of smaller aggregates}
Different drugs target aggregates of different sizes. We simulate this effect by assuming that clearance halts or reduces neurodegeneration. An enhanced clearance $\lambda_{\text{drug}}$, will specifically affect aggregates within a particular size interval, in addition to the natural (assumed to be initially subcritical) 
background clearance $\lambda_a$ of the human brain. To capture these effects in a local model representative of a single brain region \textit{in vivo} like \eqref{eq:vivo4}, we keep the clearance rates constant for most sizes but elevated dramatically across an interval $n_0 \leq i \leq n_1$: 
\begin{subequations}
\begin{align}
& \frac{{\text{d}}{{m}}}{{\text{d}}t} = 0, \\
& \frac{{\text{d}}{{p_2}}}{{\text{d}}t} = -\lambda_2 p_2 + k_{n} m^2 - 2 k_{+} m p_{2} + k_{2} \sigma(M) m^2 M, \\
& \frac{{\text{d}}{{p_i}}}{{\text{d}}t} = -\lambda_i p_i + 2 k_{+} m (p_{i-1} - p_{i}), \quad i > 2,\\
&  \lambda_i =\begin{cases}
     \lambda_a + \lambda_{\text{drug}}, & n_0 \leq i \leq n_1\\
      \lambda_a, & \text{otherwise}, \label{eq:cl_interval}
    \end{cases}
\end{align}
\end{subequations}
where 
\begin{equation}
\sigma(M) = \frac{K_M }{K_M + M^{2}}, \quad %
P=\sum_{i=2}^{\infty} p_{i}, \quad  M=\sum_{i=2}^{\infty} i p_{i},
\end{equation}
with the corresponding reduced moments system 
\begin{subequations}
\label{eq:vivo4Minterval}
\begin{align}
&{\frac{{\text{d}}{{P}}}{{\text{d}}t}}  
\,\,=  - \lambda_a P - \lambda_{\text{drug}}\sum_{i=n_1}^{n_2} p_{i} +k_{n}m_0^2+ k_{2}\, \sigma(M) m_0^2 M,
 \label{Pinvitro2bMa_interval}\\
&{\frac{{\text{d}}{{M}}}{{\text{d}}t}}  
= - \lambda_a M - \lambda_{\text{drug}}\sum_{i=n_1}^{n_2} i p_{i} + 2 k_{n}m_0^2+{2 k_{+}m_0 P} + {2}k_{2}\, \sigma(M) m_0^2 M. \label{Minvitro2bMa_interval}
\end{align}
\end{subequations}

\begin{figure}[h!]
    \centering
    \includegraphics[width = .7\textwidth, trim={1.3cm 0 0.5cm 0},clip]{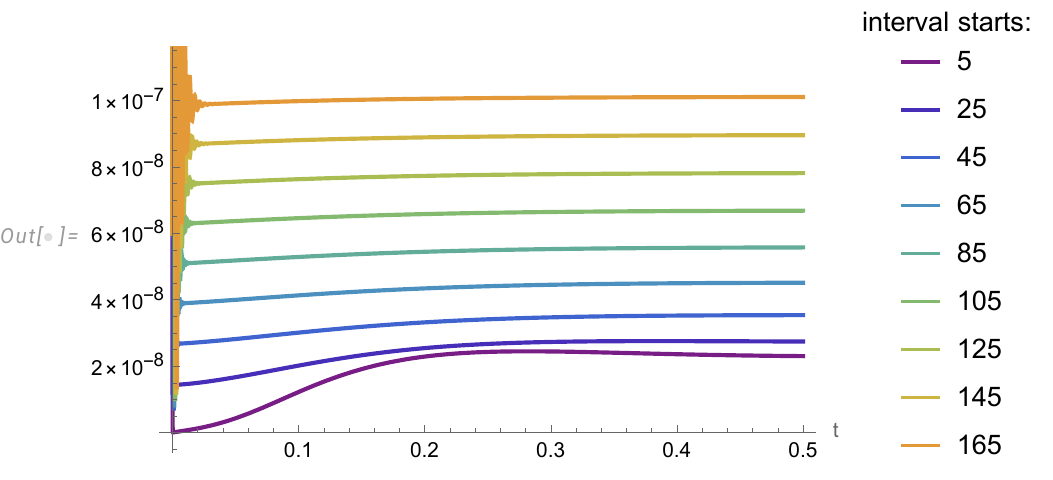}
    \caption{Total toxic mass evolution with heightened clearance as in \eqref{eq:cl_interval} with $n_0$ specified in the legend, and $n_1 = n_0+10$. The equilibrium solution for toxic mass $M_\infty$ is lower when smaller aggregates are targeted with heightened clearance. Here $\lambda_a = 10$ and $\lambda_{\text{drug}} = 10^5$ and parameters are as found in \Cref{tab:parameters}.}
    \label{fig:drug_interval}
\end{figure}

As seen in \Cref{fig:drug_interval}, we perform a sweep of intervals $(n_0,n_1)$ and observe that targeting aggregates of smaller size decreases the total toxic mass in the system the most. 
This underscores the fundamental role of secondary nucleation in forming plaques and indicates that smaller aggregates constitute the largest toxic mass in the system.

An explanation of the observed results in \Cref{fig:drug_interval} can be understood by considering the fixed point solutions. Specifically, 
a steady state $p_i^*$ for $i> 2$ satisfies the recurrence relation
\begin{equation}
    p_i^* = \delta_i p_{i-1}^*,\quad \delta_i = \frac{2k_+ m_0}{\lambda_i + 2k_+ m_0},\quad i>2,
\end{equation}
where each recursive effect is dependent on $i$. Consequently, the fixed point of the $i$-mer is 
\begin{equation}
    p_i^* = \triangle_i p_2^*, \quad \text{where } \quad 
    \triangle_i =
\begin{cases}
    \prod_{j=3}^i \delta_j & \text{for } i>2, \\
    1 & \text{for } i = 2.
\end{cases}
\end{equation}
Removing $i$-mers of sizes $n_1\leq i \leq n_2$ gives the toxic mass fixed point,
\begin{equation}
    M_* = \triangle p_2^* - \sum_{k=n_1}^{n_2} k \triangle_k p_2^*,
\end{equation}
where we define $\triangle = \sum_{k=2}^\infty k \triangle_k$, and
\begin{equation}
    p_2^*= \frac{\sqrt{23} \sqrt{-\lambda_2 - 2 k_{+} m_0 + k_2 m_0^2 \triangle}}{1000 \sqrt{
 \lambda_2 \triangle^2 + 2 k_p m_0 \triangle^2}}.
\end{equation}
It remains to prove that 
\begin{equation}
    \sum_{k=n_1}^{n_2} k \triangle_k p_2^* > \sum_{k=n_3}^{n_4} k \triangle_k p_2^*,
\end{equation}
or
\begin{equation}
    (a+\lambda n_1)\Big( \frac{a}{a + \lambda} \Big)^{n_1} - (a+\lambda n_2)\Big( \frac{a}{a + \lambda} \Big)^{n_2} > (a+\lambda n_3)\Big( \frac{a}{a + \lambda} \Big)^{n_3} - (a+\lambda n_4)\Big( \frac{a}{a + \lambda} \Big)^{n_4},
\end{equation}
for $n_1<n_2<n_3<n_4$ where $a = 2 k_+ m_0$. Expanding the brackets, to leading order 
\begin{equation}
    n_1 \xi^{n_1} - n_2 \xi^{n_2} > n_3 \xi^{n_3} - n_4 \xi^{n_4},
\end{equation}
where $\xi = a/\lambda$.
Since there is a fixed difference between intervals, say $d$, 
\begin{equation}
    n_1 \xi^{n_1} - n_2 \xi^{n_2} > \xi^{d} ((n_1+d) \xi^{n_1} - (n_2+d) \xi^{n_2} ),
\end{equation}
so we require
\begin{equation}
    (n_1 - \xi^{d} (n_1+d))\xi^{n_1}> (n_2 - \xi^{d} (n_2+d))\xi^{n_2}.
\end{equation}
Further, there is also a fixed interval size $n_2 = n_1 + n_d$ so
\begin{equation}
    (n_1 - \xi^{d} (n_1+d))\xi^{n_1}> (n_1 - \xi^{d} (n_1+d))\xi^{n_1}\xi^{n_d} + (n_d - \xi^{d} (n_d+d))\xi^{n_1}\xi^{n_d},
\end{equation}
that is, 
\begin{equation}
    \xi^{n_d} \Bigg( 1 + \frac{n_d - \xi^{d} (n_d+d)}{n_1 - \xi^{d} (n_1+d)} \Bigg) < 1,
\end{equation}
which holds true since $\xi \ll 1$. Hence, in the constant clearance regime, removing intervals of smaller aggregates reduces the total toxic mass equilibrium more than removing the same-sized intervals of larger aggregates. It remains to explore the effect of removing larger intervals of larger aggregates and a critical interval size, which would be more valuable than removing lower-weight oligomers. 

Given the pronounced toxicity of oligomeric intermediates relative to larger fibrillar species, actively eliminating oligomers has a dual impact. By targeting smaller aggregates, we significantly reduce the total toxic mass to the greatest extent possible. This approach aligns with our primary objective of addressing the most toxic species generated during aggregation. As a reference, aducanumab partially targets oligomers but mostly clears insoluble amyloid plaques \cite{tolar2020aducanumab}.

\subsection{Optimal drug administration strategies}

A variety of dosing strategies can be analysed and simulated directed in the nucleation-aggregation-clearance model with \textit{in vivo} effects \eqref{eq:vivo4}. In this section, we formulate an optimisation problem with the objective of minimising the accumulation of toxic mass in a given time period, over the parameters of the functional form of drug-induced clearance $\lambda(t)$. This allows us to determine the optimal frequency and volume of doses, while adhering to drug toxicity constraints, in the context of a computational drug trial.

Assuming a linear relationship between antibody concentration and clearance as proposed by \cite{mazer2023development}, we take the concentration of the drug administered to be directly proportional to the drug-enhanced clearance increment $\lambda_{\text{drug}}$. That is, $\lambda_{\text{drug}}(C_p) = L C_p$ where $C_p$ is drug concentration and $L$ is constant. 
In the first instance, we assume that drug administration is spaced equally and decays exponentially, being cleared from the brain as a non-aggregating particle. Emulating the drug concentration profiles described in \cite[Figure 1]{mazer2023development}, we simulate a dosing regime through the following clearance profile:
\begin{equation}
    \lambda(t) = \lambda_{\text{drug}} e^{-A\mod{(t,\,B)}} + \lambda_a,
    \label{eq:drugprofile}
\end{equation}
in \eqref{eq:vivo4}. Here $A$ is the clearance rate of the drug from the brain, $B$ is the time period (days) between drug doses, $\lambda_a$ is the background natural clearance, and $\lambda_{\text{drug}}$ is the drug-enhanced increase in clearance that is assumed to be proportional to the dose. 

\begin{figure}[h!]
  \centering
  \subfloat[Drug induced clearance profiles]{
    \includegraphics[height =.2\textheight, trim={0 0 0 0}, clip]{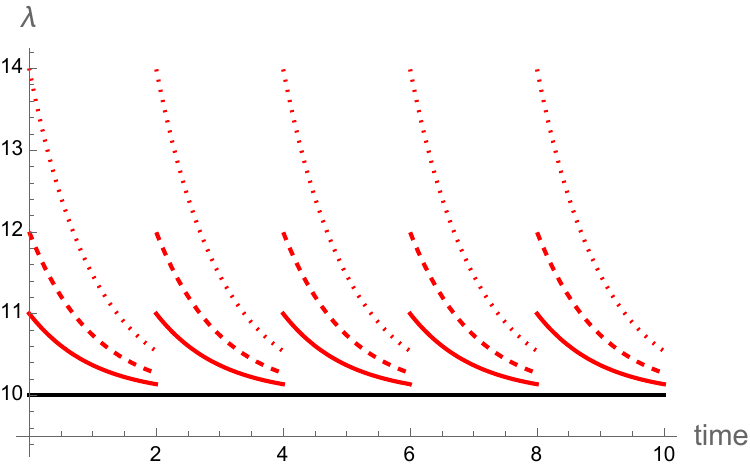}
  }
  \subfloat[Toxic mass dynamics]{
    \includegraphics[height =.2\textheight, trim={1.3cm 0 2.7cm 0}, clip]{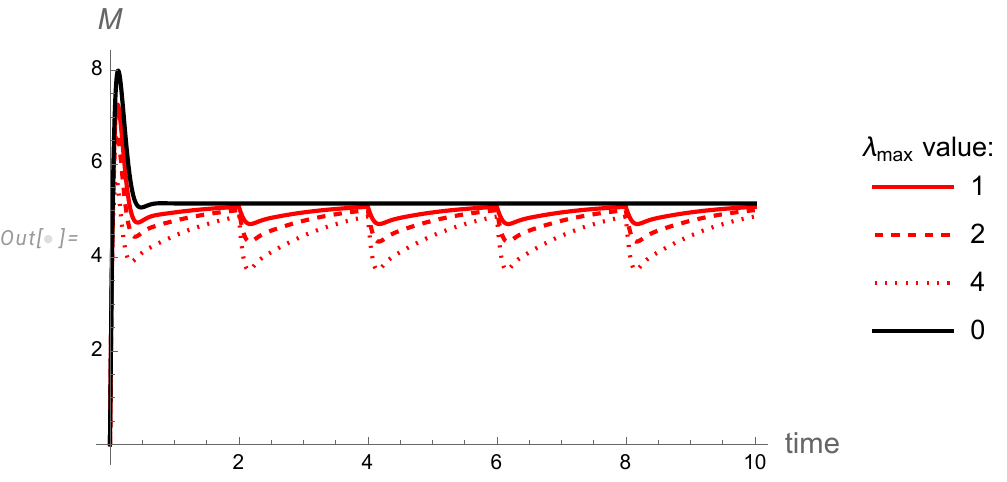}
    \includegraphics[height =.2\textheight, trim={0 0 0 0}, clip]{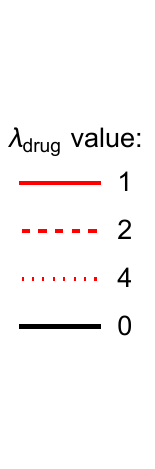}
  }
  \caption{ The drug induced clearance profiles (a) according to  \eqref{eq:drugprofile}, with $\lambda_a = 10$, $B=2$ and various $\lambda_{\text{drug}}$ as specified in the legend. The no drug case is displayed by the black curves. The corresponding toxic mass $M$ evolution shows a sinusoidal steady state in the presence of drugs. The system is \eqref{eq:vivo4M} with parameter values as in \Cref{tab:parameters}. \label{fig:drug1_profiles}}
\end{figure}

To demonstrate the effects of this dosing regime, clearance profiles \eqref{eq:drugprofile} and the response in toxic mass evolution according to \eqref{eq:vivo4} are displayed in \Cref{fig:drug1_profiles} corresponding to three dosing strategies. In this example, for illustrative purposes, the black curves show the no-treatment case, and only $\lambda_{\text{drug}}$ (the dosage) varies across strategies with fixed $A = 1$, $B=2$, and ambient clearance $\lambda_a = 10$. We observe a sinusoidal steady state in the presence of drugs, with the steady state's peak reducing and the steady state's range increasing, with $\lambda_{\text{drug}}$ increasing. 

Considering the effects of varying both $B$ and $\lambda_{\text{drug}}$ on misfolded protein mass, herein referred to as the \textit{regime parameters}, the best strategy is to maximise the frequency and volume of treatments. However, there is a cost involved in raising $\lambda_{\text{drug}}$ (dosage) and reducing $B$ (time between each dose) due to the toxicity of antibodies to the brain environment, such as a significant risk of amyloid-related imaging abnormality (ARIA) 
\cite{withington2022amyloid,cummings2021aducanumab}. With this in mind, the following optimisation problem naturally arises. 

Given a trial period of say $t=28$ days, we aim to identify the optimal number of days between drug administration $B$ and volume of drugs given that is proportional to drug-induced clearance $\lambda_{\text{drug}}$. We supplement the problem with the constraint that the total mass of drugs (mg) administered during this period invokes a limited integrated clearance increment, denoted $C_{\text{max}}$. For simplicity, we choose $A=1$, assuming that $A$ remains constant as clearance increases. We optimise $\lambda_{\text{drug}}$ and $B$ such that the average toxic mass across one cycle of drug administration, i.e. the average of the steady states shown in \Cref{fig:drug1_profiles}(b), $\Bar{M}$, is minimised. In full, the optimisation problem reads 
\begin{align}
    &\text{Minimise} \quad \Bar{M} \quad \text{over} \quad 0.2<B<t_{\text{max}}, \quad 0.2<\lambda_{\text{drug}}<80\\
    &\text{subject to} \\
     & \int_0^{t_{\text{max}}} (\lambda_{\text{drug}} e^{-A\mod{(t,\,B)}}) \, \text{d}t =  C_{\text{max}}.
    \label{constraints}
\end{align}
\begin{figure}[h!]
  \centering
  \subfloat[Average toxic mass with drugs $\Bar{M}$]{
    \includegraphics[height =0.28\textheight, trim={1.5cm 0 0 0}, clip]{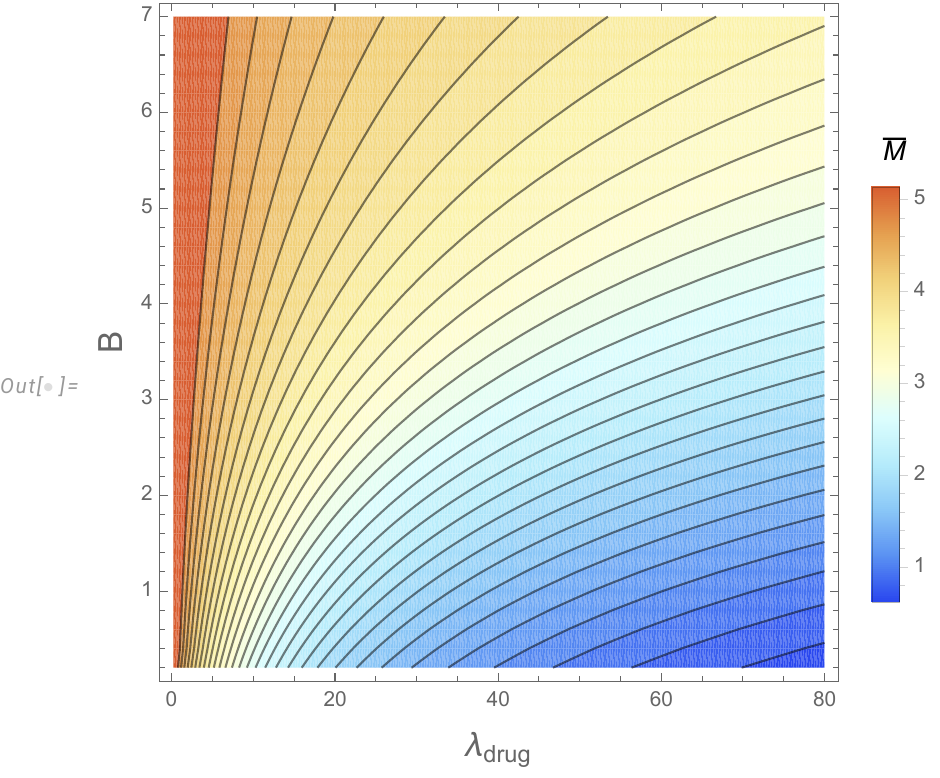}
  }
  \subfloat[Maximum drug toxicity $C_{\text{max}}$]{
    \includegraphics[height =0.28\textheight, trim={1.5cm 0 0 0}, clip]{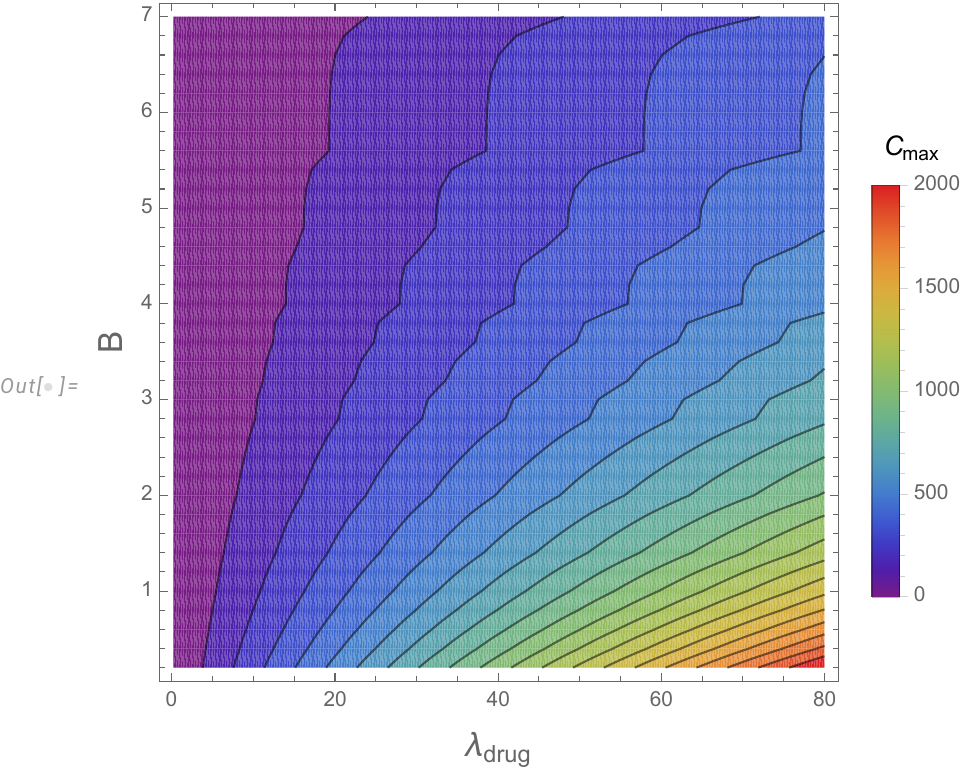}
  }
  \caption{Contour maps showing the average toxic mass $\Bar{M}$ (a) and toxic load of the drug $C_{\text{max}}$ (b) as functions of the drug dosing regime parameters in \eqref{eq:drugprofile}: the drug induced maximum clearance increment $\lambda_{\text{drug}}$, and the frequency of treatments $B$. \label{fig:optimisingdosage1}}
\end{figure}

\begin{figure}[h!]
    \centering
    \includegraphics[height = .44\textheight, trim={0 0 3cm 0},clip]{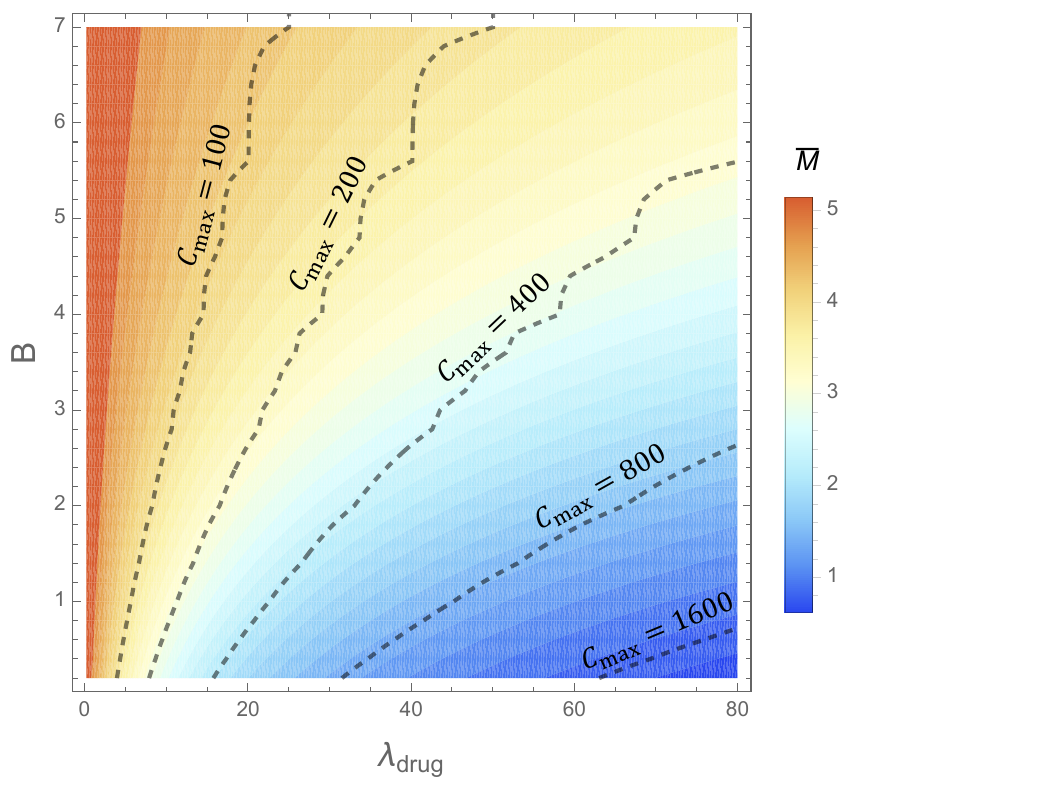}
    \caption{Contour plot with the colours displaying the average toxic mass $\Bar{M}$ accumulated in the presence of a drug with different dosing regimes defined by varying parameters $\lambda_{\text{drug}}$ correlated directly with dosage size, and $B$ the frequency of doses. Toxic mass is found numerically by solving \eqref{eq:vivo4} with \eqref{eq:drugprofile}. Solutions which satisfy the constraints outlined in \eqref{constraints}, namely with regime parameters invoking the specific toxic load $C_{\text{max}}$, for different values of $C_{\text{max}}$, lie along the contour lines.}
    \label{fig:3dlistplot}
\end{figure}

Contour plots of the average toxic mass $\Bar{M}$ for varying regime parameters are shown in \Cref{fig:optimisingdosage1}(a), with the corresponding drug toxicity $C_{\text{max}}$ shown in \Cref{fig:optimisingdosage1}(b). We identify the optimal dosing regime, subject to a specific drug toxicity constraint $C_{\text{max}}$ by consulting the combination of \Cref{fig:optimisingdosage1}(a) and (b) in \Cref{fig:3dlistplot}. Solutions satisfying the constraints $C_{\text{max}}\approx 100,\,200,400,800,1600$ in \eqref{constraints} lie along the respective dashed lines, and the colour function displays the relative average toxic mass equilibrium $\Bar{M}$. Subject to these constraints, solutions with the lowest average toxic mass $\Bar{M}$ correspond to the lowest spacing between doses, $B\to0$. 

The optimal strategy, therefore, is the trivial one: to take $B\to 0$, corresponding to a constant supply of drugs. Of course, this could be impractical or detrimental to a patient's quality of life. Given the potential impracticality of constant drug supply, the key question is at what cost, in terms of toxic mass $\Bar{M}$, can drug doses be distributed (increasing $B$). For reference, consider regime parameters corresponding to $C_{\text{max}} = 100$ in \eqref{constraints}. The optimal solution, a constant supply of drugs with $\lambda_{\text{drug}}= 100/28$, results in the minimum allowable average toxic mass of $\Bar{M}=3.79$. Other solutions satisfying the same drug toxicity constraint, but not necessarily minimising $\Bar{M}$, are shown in \Cref{fig:3dlistplot} and demonstrate some flexibility in dosing strategies. For example, $\lambda_{\text{drug}} = 5.6$ once a day results in the only slightly increased $\Bar{M}=3.81$. A solution with seven days between doses of $\lambda_{\text{drug}} = 25$ results in a similar average toxic steady state of $\Bar{M}=4.22$, albeit with a higher variance around this equilibrium. This implies some leniency on dosing frequency. Overall, as seen in \Cref{fig:3dlistplot}, $\Bar{M}$ does not greatly vary along the constant $C_{\text{max}}$ contours; the contours of $\Bar{M}$ and $C_{\text{max}}$ in \Cref{fig:optimisingdosage1}(a) and (b) follow similar trends. The amount that $\Bar{M}$ varies with $B$ with $C_{\text{max}}$ constant is dependent on $C_{\text{max}}$, as observed by the higher range in $\Bar{M}$ values along the $C_{\text{max}}=400$ contour compared to the   $C_{\text{max}}=100$ contour. Therefore, depending on $C_{\text{max}}$ and the desired reduction in toxic mass, increasing the days between doses can lead to insignificant increases in toxic mass, with potentially significant improvements in quality of life. 
\begin{figure}[h!]
    \centering
    \includegraphics[width=.5\textwidth, trim={0 0 0 0},clip]{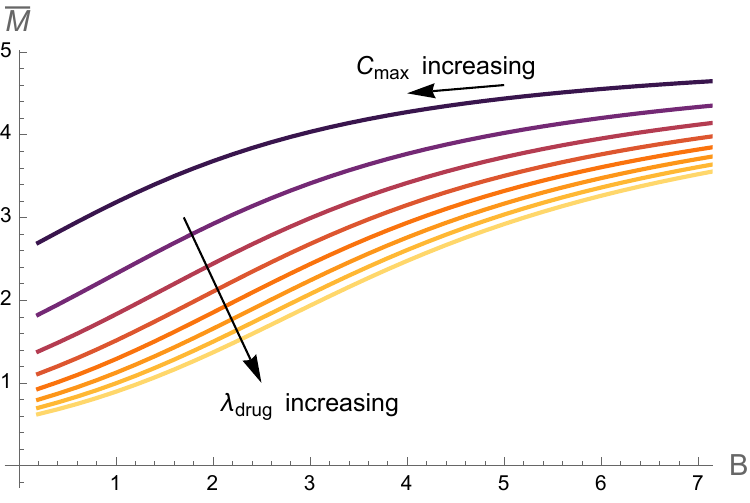}
    \caption{The evolution of the average toxic mass $\Bar{M}$ with constant drug dosage corresponding to a clearance increment of $\lambda_{\text{drug}}$ and varying frequencies of administration $B$. Constant $\lambda_{\text{drug}}$ values range from 10 to 80 with differences of 10 between curves.
    \label{fig:varyingB}}
\end{figure}

In addition, as shown by \Cref{fig:varyingB}, keeping $\lambda_{\text{max}}$ constant and lowering $B$ can have a high impact on reducing $\Bar{M}$ for smaller values of $B$ but is less impactful for higher $B$. Further, increasing $\lambda_{\text{drug}}$ is less impactful for higher values of $\lambda_{\text{drug}}$ (shown by the lighter curves). Together, these results hold potential in determining dosing strategies in personalised medicine, with initial conditions of \eqref{eq:vivo4} reflecting the initial toxic loads of patients and parameters in \eqref{eq:drugprofile} tunable to the desired reduction in toxic load. 

\section{Discussion}

Theoretical research into toxic mass accumulation in neurodegenerative diseases has thus far exclusively consisted of either detailed \textit{in vitro} analysis of aggregation kinetics or \textit{in vivo} studies on the effects of transport and clearance in macroscale models on networks, with the former being validated against experimental observations and the latter with structural data. Both aspects have proved essential for our understanding of AD pathology. However, there is a cloudy middle ground restricting our understanding of aggregation inside the complex domain of the human brain. The development of a nucleation-aggregation model which includes \textit{in vivo} effects such as protein production, saturation, transport and, most importantly, the clearance of solutes, remains a task that has only been addressed mechanistically in two prior studies \cite{fornari2020spatially, meislclearance2020}. Fornari \textit{et al.} (2019) \cite{fornari2020spatially} used a Smoluchowski-type model to analyse the effect of aggregate transport at the brain scale, while Thompson \textit{et al.} (2021) \cite{meislclearance2020} considered the effects of a bulk constant clearance at the microscopic scale. Our work combines both approaches within the same mathematical framework, establishing a general picture for the self-assembly of A$\beta$ in the aggregation cascade \textit{in vivo} across the connectome. We have derived a class of nucleation-aggregation-clearance models for the homogeneous dynamics and spatial progression of A$\beta$ in the brain, tracking the evolution of different size aggregates with distinct properties instrumental in neurodegeneration. This new class of models provides a therapeutic modelling platform reflective of the human brain for the simulation of potential therapies at the spatial and temporal scales of the disease. Notably, we reveal the therapeutic effect of enhanced clearance behind the success of recent drug trials. The dependence of drug efficacy on initial clearance, along with the characterisation of the drug's effects purely in terms of clearance, further emphasises the pivotal role of clearance in neurodegeneration.

\subsection{Brain scale models of nucleation, aggregation, clearance and transport}

The first primary aim of this work was to develop a comprehensive and entirely mechanistic model encompassing the aforementioned \textit{in vivo} effects, along with physiological effects not previously modeled, with a 
focus on quantifying the role of dynamic clearance of aggregates in neurodegeneration at the brain scale. To address the complexities of developing a model representative of a single brain region, we advance current state-of-the-art chemical kinetic models by incorporating physiological effects observed in experiments. Then, we include strong transport anisotropy in the brain through the graph Laplacian to simulate aggregate dynamics across the spectrum of aggregate sizes at the brain scale. 
The study of such systems is guided by the homogeneous case for which both total mass evolution and aggregates' size distribution are obtained analytically, and bifurcation analysis reveals critical clearance regimes dependent on kinetic rate constants.

At the local ROI level, analysis of the homogeneous nucleation-aggregation-clearance models \eqref{eq:vivo4}
reveals that the clearance rate dictates the timescale of the nucleated transition to an aggregated brain region and the final toxic mass equilibrium. In the absence of a simulated seed, an increase in toxic mass is characterised by an early phase of seeding depending on primary nucleation, followed by a period of linear growth mostly controlled by aggregation of monomers onto the fibril, and then a phase of autocatalytic growth dominated by the multiplicative processes of secondary nucleation and elongation. Once enough toxic mass accumulates, clearance dominates, and the toxic mass approaches a saturated steady state. Importantly the critical clearance rate, above which all aggregates clear from a brain region, holds significant therapeutic potential. Recovering similar results to \cite{meislclearance2020}, the critical clearance formulae depend upon the assumed clearance form relative to aggregate size, with different dependencies on the aggregation rate constants. When analysing the local system with dynamic clearance given by \eqref{eq:vivo6}, local clearance reduces to the basal clearance capacity at rates relative to the difference between initial clearance and basal clearance. After this point the dynamics of the system are equivalent to those of the homogeneous nucleation-aggregation-clearance models \eqref{eq:vivo4} so we recover the same critical \textit{basal} clearance regimes as for clearance in \eqref{eq:vivo4}.

To model the interplay between aggregation and clearance at the brain-scale, the local dynamics are coupled with anisotropic transport across the connectome 
in \eqref{super_general_system}. As identified first by \cite{fornari2020spatially}, the progression of the disease consists of an initial stage that develops at the seeded node, followed by primary infection in connected nodes, mainly driven by diffusion rather than aggregation kinetics. These nodes distribute seeds, and due to the brain's small-world network structure, secondary infection develops in most nodes soon after primary infection. In the growth-dominated regime representative of the toxic protein cascade in a human brain, local dynamics dominate once a node receives a seed by primary or secondary infection. We thus find that the analytical results of the homogeneous model \eqref{eq:vivo4} are a good approximation for the connected dynamics at the brain scale. Furthermore, the brain-scale dynamics appear invariant to the diffusion of aggregates larger than nuclei, indicating that the transport of dimers is the primary factor driving pathology, and size-dependent diffusion does not significantly alter the dynamics.

\subsection{Therapeutic modelling insights}

The pivotal field of Alzheimer's drug design commands extensive research activity. Yet, drug trials and drug performance projections rely on data without providing comprehensive insights into the mechanisms of drug action. Current microscopic chemical kinetic models shed light on how drugs influence the aggregation process \cite{linse2020kinetic}, while empirical compartmental models project the longer-term effects of anti-amyloid treatments by fitting data from clinical trials to semi-mechanistic macroscale models \cite{mazer2023development}, with no consideration of the drug action at the microscale. Our \textit{in vivo} nucleation-aggregation-clearance model aims to bridge the gap between these two approaches, extending \textit{in vitro} observations to decades of AD pathology at the brain scale, considering both the increase in clearance due to drugs and the inhibition of specific steps in the aggregation process. We thus present the first fully mechanistic model of brain-scale \textit{in vivo} aggregation growth, clearance and transport that can describe the effects of therapeutic agents both analytically and through direct simulation, aid drug performance projection, and lead exploratory theoretical studies to ultimately guide drug design.
We thus quantified the action of potential drugs on protein aggregation \textit{in vivo}, and conducted a preliminary study on the most favorable treatment strategies, concerning optimal aggregate sizes to target, specific microscopic steps of aggregation to inhibit, and the frequency and volume of antibody dosage.

Specifically, we modeled the mechanism of action of treatments in the homogeneous system \eqref{eq:vivo4_adu} and its corresponding brain-scale extensions to see the effect of the antibodies on the disease dynamics. For aducanumab, we modeled the pharmacodynamic property of a decreased secondary nucleation rate constant due to the drug binding to fibrils, thus inhibiting the aggregation process \cite{linse2020kinetic} and mediating enhanced clearance, as considered by \cite{mazer2023development}. First focusing solely on the kinetic drug effects we observe a consistent pattern across all models: a decrease in critical clearance, leading to an improvement in the effectiveness of compromised local clearance. This enables a patient under treatment to completely transition to a healthy state or experience a lessened toxic load, depending on the initial clearance at the time of commencing therapy.
Mathematically, the inhibitory actions of potential therapeutic agents on the aggregation process transforms the aggregated equilibrium (\Cref{fig:toxic_equi}). Thus through steady-state analysis, we identify that the primary therapeutic effect of the drugs is a reduction in critical clearance, increasing the system's overall \textit{effective clearance}.

Considering the drug's effect in the microscopic aggregation model \eqref{eq:vivo4_adu} coupled to transport allows us to track the response in the evolution of distinct aggregate sizes across the connectome. The effects of the kinetic action of treatment are the same: a transformation of the toxic steady-state solution in the presence of drugs (\Cref{fig:toxic_equi}), with the potential to initiate a transition from an unhealthy patient state ($\lambda <\lambda_{\text{crit}}$) to a healthy state ($\lambda >\tilde{\lambda}_{\text{crit}}$) in which the otherwise compromised clearance rates are sufficient to completely remove toxic mass in the presence of drugs. 
The complete inhibition of aggregation and propagation shown in \Cref{fig:Abetadrugvsnodrug} is due only to the kinetic action of the drug reducing $\lambda_{\text{crit}}$ increasing the effectiveness of the brains damaged clearance mechanisms. This transition from an unhealthy state to a healthy equilibrium is dependent upon the brain's initial clearance upon drug administration, underscoring the need for treatments early on in disease progression when $\tilde{\lambda}_{\text{crit}}< \lambda<\lambda_{\text{crit}}$. 
The impact of drugs such as aducanumab on secondary nucleation thus has a profound effect, employing a decrease in the aggregation kinetic rate constants thereby improving the effectiveness of brain clearance. Our most critical finding is that the drug's kinetic effect alone can be sufficient for a transition from a brain with detrimental levels of toxic mass to a healthy state. In addition, considering enhanced clearance due to drug binding in combination with the kinetic effect of the drugs leads to a more pronounced reduction in long-term toxic mass accumulation. Significantly, in the presence of drugs, the critical toxic seed as discussed in \cite{brennanpreprint} is increased, so it is less likely that a dynamic clearance will drop below critical.

Ultimately, our brain-scale nucleation-aggregation-clearance models open up avenues for the mathematically-informed therapeutic strategies to control pathological protein aggregation, addressing questions that current modelling studies fail to address at the brain scale. We have identified that a reduction in $\lambda_{\text{crit}}$ is successful in curing patients. Quantifying the actions of drugs in this way allows us to identify which rate constants might inhibit aggregation to provide the largest reduction in  $\lambda_{\text{crit}}$. For example, if considering the size dependant clearance $\lambda_i = \lambda_0/i$ to be reflective of the human body, the form of the critical clearance rate \eqref{eq:lcrit_decreasing_cl} suggests that targeting the elongation processes would be most effective in reducing the critical clearance and most vitally translating the fixed point solution curve.  
Further, we performed a sweep of an interval of increased clearance and uncovered that targeting smaller aggregates is more impactful than targeting larger ones. Additionally, we studied an optimization problem to demonstrate how the models can inform dosing regimens, as described and simulated in \cite{mazer2023development}. The results suggest that more frequent dosing is preferable for reducing average toxic mass and variance, although combinations of higher, less frequent doses can achieve the same steady states. A more in-depth analysis of the application of these models for pharmacological use is of great interest for future studies.

\section*{Acknowledgements}
The work of A.G. was supported by the Engineering and
Physical Sciences Research Council under Research grant EP/R020205/1. This publication is based on work supported by the EPSRC Centre For Doctoral Training in Industrially Focused Mathematical Modelling (EP/L015803/1) in collaboration with Simula Research Laboratory.
For the purpose of Open Access, the authors will apply a CC BY public copyright license to any Author Accepted Manuscript (AAM) version arising from this submission. 
The help of Travis Thompson and Hadrien Oliveri with numerical and conceptual issues is gratefully acknowledged.

\section*{Data Availability}
The manuscript has no associated data.


\begin{thebibliography}{10}

    \bibitem{fowler2006functional}
    D.~M. Fowler, A.~V. Koulov, C.~Alory-Jost, M.~S. Marks, W.~E. Balch, and J.~W. Kelly.
    \newblock Functional amyloid formation within mammalian tissue.
    \newblock {\em PLoS Biology}, 4(1):e6, 2006.
    
    \bibitem{maji2009functional}
    S.~K. Maji, M.~H. Perrin, M.~R. Sawaya, S.~Jessberger, K.~Vadodaria, R.~A. Rissman, P.~S. Singru, K.~P.~R. Nilsson, R.~Simon, D.~Schubert, et~al.
    \newblock Functional amyloids as natural storage of peptide hormones in pituitary secretory granules.
    \newblock {\em Science}, 325(5938):328--332, 2009.
    
    \bibitem{bleem2017structural}
    A.~Bleem and V.~Daggett.
    \newblock Structural and functional diversity among amyloid proteins: agents of disease, building blocks of biology, and implications for molecular engineering.
    \newblock {\em Biotechnology and Bioengineering}, 114(1):7--20, 2017.
    
    \bibitem{knowles2016amyloid}
    T.~P.~J. Knowles and R.~Mezzenga.
    \newblock Amyloid fibrils as building blocks for natural and artificial functional materials.
    \newblock {\em Advanced Materials}, 28(31):6546--6561, 2016.
    
    \bibitem{chiti2006protein}
    F.~Chiti and C.~M. Dobson.
    \newblock Protein misfolding, functional amyloid, and human disease.
    \newblock {\em Annual Review of Biochemistry}, 75:333--366, 2006.
    
    \bibitem{chiti2017protein}
    F.~Chiti and C.~M. Dobson.
    \newblock Protein misfolding, amyloid formation, and human disease: a summary of progress over the last decade.
    \newblock {\em Annual Review of Biochemistry}, 86:27--68, 2017.
    
    \bibitem{knowles2014amyloid}
    T.~P.~J. Knowles, M.~Vendruscolo, and C.~M. Dobson.
    \newblock The amyloid state and its association with protein misfolding diseases.
    \newblock {\em Nature Reviews Molecular Cell Biology}, 15(6):384--396, 2014.
    
    \bibitem{wortmann2012dementia}
    M.~Wortmann.
    \newblock Dementia: a global health priority-highlights from an {ADI} and {W}orld {H}ealth {O}rganization report.
    \newblock {\em Alzheimer's Research \& Therapy}, 4:1--3, 2012.
    
    \bibitem{lansbury1996reductionist}
    P.~T. Lansbury.
    \newblock A reductionist view of {A}lzheimer's disease.
    \newblock {\em Accounts of Chemical Research}, 29(7):317--321, 1996.
    
    \bibitem{goedert2017like}
    M.~Goedert, M.~Masuda-Suzukake, and B.~Falcon.
    \newblock Like prions: the propagation of aggregated tau and $\alpha$-synuclein in neurodegeneration.
    \newblock {\em Brain}, 140(2):266--278, 2017.
    
    \bibitem{frost2010prion}
    B.~Frost and M.~I. Diamond.
    \newblock Prion-like mechanisms in neurodegenerative diseases.
    \newblock {\em Nature Reviews Neuroscience}, 11(3):155--159, 2010.
    
    \bibitem{jucker2018propagation}
    M.~Jucker and L.~C. Walker.
    \newblock Propagation and spread of pathogenic protein assemblies in neurodegenerative diseases.
    \newblock {\em Nature Neuroscience}, 21(10):1341--1349, 2018.
    
    \bibitem{olsson2018prion}
    T.~T. Olsson, O.~Klementieva, and G.~K. Gouras.
    \newblock Prion-like seeding and nucleation of intracellular amyloid-\textbeta.
    \newblock {\em Neurobiology of Disease}, 113:1--10, 2018.
    
    \bibitem{goedert2015alzheimer}
    M.~Goedert.
    \newblock {{A}lzheimer's and {P}arkinson's diseases: The prion concept in relation to assembled {A}\textbeta{}, tau, and \textalpha{}-synuclein}.
    \newblock {\em Science}, 349(6248):1255555, 2015.
    
    \bibitem{mudher2017evidence}
    A.~Mudher, M.~Colin, S.~Dujardin, M.~Medina, I.~Dewachter, S.~M.~A. Naini, E.-M. Mandelkow, E.~Mandelkow, L.~Bu{\'e}e, M.~Goedert, et~al.
    \newblock What is the evidence that tau pathology spreads through prion-like propagation?
    \newblock {\em Acta Neuropathologica Communications}, 5(1):99, 2017.
    
    \bibitem{prusiner1998prions}
    S.~B. Prusiner.
    \newblock Prions.
    \newblock {\em Proceedings of the National Academy of Sciences}, 95(23):13363--13383, 1998.
    
    \bibitem{hardy1992alzheimer}
    J.~A. Hardy and G.~A. Higgins.
    \newblock {A}lzheimer's disease: the amyloid cascade hypothesis.
    \newblock {\em Science}, 256(5054):184--186, 1992.
    
    \bibitem{hardy1991amyloid}
    J.~Hardy and D.~Allsop.
    \newblock Amyloid deposition as the central event in the aetiology of {A}lzheimer's disease.
    \newblock {\em Trends in Pharmacological Sciences}, 12:383--388, 1991.
    
    \bibitem{selkoe2016amyloid}
    D.~J. Selkoe and J.~Hardy.
    \newblock The amyloid hypothesis of {A}lzheimer's disease at 25 years.
    \newblock {\em EMBO Molecular Medicine}, 8(6):595--608, 2016.
    
    \bibitem{selkoe2001alzheimer}
    D.~J. Selkoe.
    \newblock Alzheimer's disease: genes, proteins, and therapy.
    \newblock {\em Physiological Reviews}, 2001.
    
    \bibitem{hard2012inhibition}
    T.~H{\"a}rd and C.~Lendel.
    \newblock Inhibition of amyloid formation.
    \newblock {\em Journal of Molecular Biology}, 421(4-5):441--465, 2012.
    
    \bibitem{karran2022amyloid}
    E.~Karran and B.~De~Strooper.
    \newblock The amyloid hypothesis in {A}lzheimer disease: new insights from new therapeutics.
    \newblock {\em Nature Reviews Drug Discovery}, 21(4):306--318, 2022.
    
    \bibitem{hong2018diffusible}
    W.~Hong, Z.~Wang, W.~Liu, T.~T. O’Malley, M.~Jin, M.~Willem, C.~Haass, M.~P. Frosch, and D.~M. Walsh.
    \newblock Diffusible, highly bioactive oligomers represent a critical minority of soluble {A}$\beta$ in {A}lzheimer’s disease brain.
    \newblock {\em Acta Neuropathologica}, 136:19--40, 2018.
    
    \bibitem{walsh2020amyloid}
    D.~M. Walsh and D.~J. Selkoe.
    \newblock Amyloid $\beta$-protein and beyond: The path forward in {A}lzheimer’s disease.
    \newblock {\em Current Opinion in Neurobiology}, 61:116--124, 2020.
    
    \bibitem{linse2020kinetic}
    S.~Linse, T.~Scheidt, K.~Bernfur, M.~Vendruscolo, C.~M. Dobson, S.~I.~A. Cohen, E.~Sileikis, M.~Lundqvist, F.~Qian, T.~O’Malley, et~al.
    \newblock Kinetic fingerprints differentiate the mechanisms of action of anti-{A}$\beta$ antibodies.
    \newblock {\em Nature Structural \& Molecular Biology}, 27(12):1125--1133, 2020.
    
    \bibitem{rumble1989amyloid}
    B.~Rumble, R.~Retallack, C.~Hilbich, G.~Simms, G.~Multhaup, R.~Martins, A.~Hockey, P.~Montgomery, K.~Beyreuther, and C.~L. Masters.
    \newblock Amyloid {A}4 protein and its precursor in {D}own's syndrome and {A}lzheimer's disease.
    \newblock {\em The New England Journal of Medicine}, 320(22):1446--1452, 1989.
    
    \bibitem{bacyinski2017paravascular}
    A.~Bacyinski, M.~Xu, W.~Wang, and J.~Hu.
    \newblock The paravascular pathway for brain waste clearance: current understanding, significance and controversy.
    \newblock {\em Frontiers in Neuroanatomy}, 11:101, 2017.
    
    \bibitem{tarasoff2015clearance}
    J.~M. Tarasoff-Conway, R.~O. Carare, R.~S. Osorio, L.~Glodzik, T.~Butler, E.~Fieremans, L.~Axel, H.~Rusinek, C.~Nicholson, B.~V. Zlokovic, et~al.
    \newblock Clearance systems in the brain - implications for {A}lzheimer disease.
    \newblock {\em Nature Reviews Neurology}, 11(8):457--470, 2015.
    
    \bibitem{xin2018clearance}
    S.-H. Xin, L.~Tan, X.~Cao, J.-T. Yu, and L.~Tan.
    \newblock Clearance of amyloid beta and tau in {A}lzheimer’s disease: from mechanisms to therapy.
    \newblock {\em Neurotoxicity Research}, 34:733--748, 2018.
    
    \bibitem{dobson2013story}
    M.~Dobson.
    \newblock {\em The Story of Medicine}.
    \newblock Quercus Books, 2013.
    
    \bibitem{klein2019gantenerumab}
    G.~Klein, P.~Delmar, N.~Voyle, S.~Rehal, C.~Hofmann, D.~Abi-Saab, M.~Andjelkovic, S.~Ristic, G.~Wang, R.~Bateman, et~al.
    \newblock Gantenerumab reduces amyloid-$\beta$ plaques in patients with prodromal to moderate {A}lzheimer’s disease: a {PET} substudy interim analysis.
    \newblock {\em Alzheimer's Research \& Therapy}, 11(1):1--12, 2019.
    
    \bibitem{ten2018secondary}
    M.~Ten~Kate, S.~Ingala, A.~J. Schwarz, N.~C. Fox, G.~Ch{\'e}telat, B.~N.~M. van Berckel, M.~Ewers, C.~Foley, J.~D. Gispert, D.~Hill, et~al.
    \newblock Secondary prevention of {A}lzheimer’s dementia: neuroimaging contributions.
    \newblock {\em Alzheimer's Research \& Therapy}, 10(1):1--21, 2018.
    
    \bibitem{walker2015neurodegenerative}
    L.~C. Walker and M.~Jucker.
    \newblock Neurodegenerative diseases: expanding the prion concept.
    \newblock {\em Annual Review of Neuroscience}, 38:87--103, 2015.
    
    \bibitem{fornari2019prion}
    S.~Fornari, A.~Sch{\"a}fer, M.~Jucker, A.~Goriely, and E.~Kuhl.
    \newblock Prion-like spreading of {A}lzheimer's disease within the brain's connectome.
    \newblock {\em Journal of the Royal Society Interface}, 16(159):20190356, 2019.
    
    \bibitem{fornari2020spatially}
    S.~Fornari, A.~Sch{\"a}fer, E.~Kuhl, and A.~Goriely.
    \newblock Spatially-extended nucleation-aggregation-fragmentation models for the dynamics of prion-like neurodegenerative protein-spreading in the brain and its connectome.
    \newblock {\em Journal of Theoretical Biology}, 486:110102, 2020.
    
    \bibitem{brennanpreprint}
    G.~S. Brennan, T.~B. Thompson, H.~Oliveri, M.~E. Rognes, and A.~Goriely.
    \newblock The role of clearance in neurodegenerative diseases.
    \newblock {\em SIAM Journal on Applied Mathematics}, pages S172--S198, 2023.
    
    \bibitem{mazer2023development}
    N.~A. Mazer, C.~Hofmann, D.~Lott, R.~Gieschke, G.~Klein, F.~Boess, H.~P. Grimm, G.~A. Kerchner, M.~Baudler-Klein, J.~Smith, et~al.
    \newblock Development of a quantitative semi-mechanistic model of {A}lzheimer's disease based on the amyloid/tau/neurodegeneration framework (the {Q-ATN} model).
    \newblock {\em Alzheimer's \& Dementia}, 19(6):2287--2297, 2023.
    
    \bibitem{meisl2017scaling}
    G.~Meisl, L.~Rajah, S.~A.~I. Cohen, M.~Pfammatter, A.~{\v{S}}ari{\'c}, E.~Hellstrand, A.~K. Buell, A.~Aguzzi, S.~Linse, M.~Vendruscolo, et~al.
    \newblock Scaling behaviour and rate-determining steps in filamentous self-assembly.
    \newblock {\em Chemical Science}, 8(10):7087--7097, 2017.
    
    \bibitem{frankel2019autocatalytic}
    R.~Frankel, M.~T{\"o}rnquist, G.~Meisl, O.~Hansson, U.~Andreasson, H.~Zetterberg, K.~Blennow, B.~Frohm, T.~Cedervall, T.~P.~J. Knowles, et~al.
    \newblock Autocatalytic amplification of {A}lzheimer-associated {A}$\beta$42 peptide aggregation in human cerebrospinal fluid.
    \newblock {\em Communications Biology}, 2(1):365, 2019.
    
    \bibitem{kundel2018measurement}
    F.~Kundel, L.~Hong, B.~Falcon, W.~A. McEwan, T.~C.~T. Michaels, G.~Meisl, N.~Esteras, A.~Y. Abramov, T.~J.~P. Knowles, M.~Goedert, et~al.
    \newblock Measurement of tau filament fragmentation provides insights into prion-like spreading.
    \newblock {\em ACS Chemical Neuroscience}, 9(6):1276--1282, 2018.
    
    \bibitem{cohen2011nucleated1}
    S.~I.~A. Cohen, M.~Vendruscolo, M.~E. Welland, C.~M. Dobson, E.~M. Terentjev, and T.~P.~J. Knowles.
    \newblock Nucleated polymerization with secondary pathways. {I}. {T}ime evolution of the principal moments.
    \newblock {\em The Journal of Chemical Physics}, 135(6), 2011.
    
    \bibitem{cohen2011nucleated2}
    S.~I.~A. Cohen, M.~Vendruscolo, C.~M. Dobson, and T.~P.~J. Knowles.
    \newblock Nucleated polymerization with secondary pathways. {II}. {D}etermination of self-consistent solutions to growth processes described by non-linear master equations.
    \newblock {\em The Journal of Chemical Physics}, 135(6), 2011.
    
    \bibitem{cohen2011nucleated3}
    S.~I.~A. Cohen, M.~Vendruscolo, C.~M. Dobson, and T.~P.~J. Knowles.
    \newblock Nucleated polymerization with secondary pathways. {III}. {E}quilibrium behavior and oligomer populations.
    \newblock {\em The Journal of Chemical Physics}, 135(6), 2011.
    
    \bibitem{meisl2016quantitative}
    G.~Meisl, X.~Yang, B.~Frohm, T.~P.~J. Knowles, and S.~Linse.
    \newblock Quantitative analysis of intrinsic and extrinsic factors in the aggregation mechanism of alzheimer-associated {A}$\beta$-peptide.
    \newblock {\em Scientific Reports}, 6(1):18728, 2016.
    
    \bibitem{yang2018role}
    X.~Yang, G.~Meisl, B.~Frohm, E.~Thulin, T.~P.~J. Knowles, and S.~Linse.
    \newblock On the role of sidechain size and charge in the aggregation of a $\beta$ 42 with familial mutations.
    \newblock {\em Proceedings of the National Academy of Sciences}, 115(26):E5849--E5858, 2018.
    
    \bibitem{nilsson2002low}
    M.~R. Nilsson, M.~Driscoll, and D.~P. Raleigh.
    \newblock Low levels of asparagine deamidation can have a dramatic effect on aggregation of amyloidogenic peptides: implications for the study of amyloid formation.
    \newblock {\em Protein Science}, 11(2):342--349, 2002.
    
    \bibitem{jucker2011pathogenic}
    M.~Jucker and L.~C. Walker.
    \newblock Pathogenic protein seeding in {A}lzheimer disease and other neurodegenerative disorders.
    \newblock {\em Annals of Neurology}, 70(4):532--540, 2011.
    
    \bibitem{jucker2013self}
    M.~Jucker and L.~C. Walker.
    \newblock Self-propagation of pathogenic protein aggregates in neurodegenerative diseases.
    \newblock {\em Nature}, 501(7465):45--51, 2013.
    
    \bibitem{meislclearance2020}
    T.~B. Thompson, G.~Meisl, T.~P.~J. Knowles, and A.~Goriely.
    \newblock The role of clearance mechanisms in the kinetics of pathological protein aggregation involved in neurodegenerative diseases.
    \newblock {\em Journal of Chemical Physics}, 154(12):125101, 2021.
    
    \bibitem{meisl2014differences}
    G.~Meisl, X.~Yang, E.~Hellstrand, B.~Frohm, J.~B. Kirkegaard, S.~I.~A. Cohen, C.~M. Dobson, S.~Linse, and T.~P.~J. Knowles.
    \newblock Differences in nucleation behavior underlie the contrasting aggregation kinetics of the {A}$\beta$40 and {A}$\beta$42 peptides.
    \newblock {\em Proceedings of the National Academy of Sciences U.S.A.}, 111(26):9384--9389, 2014.
    
    \bibitem{cohen2013proliferation}
    S.~I.~A. Cohen, S.~Linse, L.~M. Luheshi, E.~Hellstrand, D.~A. White, L.~Rajah, D.~E. Otzen, M.~Vendruscolo, C.~M. Dobson, and T.~P.~J. Knowles.
    \newblock Proliferation of amyloid-$\beta$42 aggregates occurs through a secondary nucleation mechanism.
    \newblock {\em Proceedings of the National Academy of Sciences}, 110(24):9758--9763, 2013.
    
    \bibitem{oosawa1975thermodynamics}
    F.~Oosawa and S.~Asakura.
    \newblock Thermodynamics of the Polymerization of Protein.
    \newblock {\em Academic Press}, 1975.
    
    \bibitem{oosawa1970size}
    F.~Oosawa.
    \newblock Size distribution of protein polymers.
    \newblock {\em Journal of Theoretical Biology}, 27(1):69--86, 1970.
    
    \bibitem{bennett2018tau}
    R.~E. Bennett, A.~B. Robbins, M.~Hu, X.~Cao, R.~A. Betensky, T.~Clark, S.~Das, and B.~T. Hyman.
    \newblock Tau induces blood vessel abnormalities and angiogenesis-related gene expression in {P}301{L} transgenic mice and human {A}lzheimer’s disease.
    \newblock {\em Proceedings of the National Academy of Sciences}, 115(6):E1289--E1298, 2018.
    
    \bibitem{carrillo2014}
    P.~Carrillo-Mora, R.~Luna, and L.~Colin-Barenque.
    \newblock {A}myloid beta: multiple mechanisms of toxicity and only some protective effects?
    \newblock {\em Oxidative Medicine and Cellular Longevity}, 2014:795375, 2014.
    
    \bibitem{canobbio2015role}
    I.~Canobbio, A.~A. Abubaker, C.~Visconte, M.~Torti, and G.~Pula.
    \newblock Role of amyloid peptides in vascular dysfunction and platelet dysregulation in {A}lzheimer’s disease.
    \newblock {\em Frontiers in Cellular Neuroscience}, 9:65, 2015.
    
    \bibitem{michalicova2020tau}
    A.~Michalicova, P.~Majerova, and A.~Kovac.
    \newblock Tau protein and its role in blood--brain barrier dysfunction.
    \newblock {\em Frontiers in Molecular Neuroscience}, 13:570045, 2020.
    
    \bibitem{brennan_goriely_chapter}
    G.~S. Brennan and A.~Goriely.
    \newblock An introduction to network models of neurodegenerative diseases.
    \newblock In J.~Dokken, K.-A. Mardal, M.~E. Rognes, L.~M. Valnes, and V.~Vinje, editors, {\em Mathematical modeling of the human brain: from deep learning to glymphatics}. Springer, 2025.
    \newblock (accepted).
    
    \bibitem{bertsch2017alzheimer}
    M.~Bertsch, B.~Franchi, N.~Marcello, M.~C. Tesi, and A.~Tosin.
    \newblock Alzheimer's disease: a mathematical model for onset and progression.
    \newblock {\em Mathematical Medicine and Biology: A Journal of the IMA}, 34(2):193--214, 2017.
    
    \bibitem{putra2023front}
    P.~Putra, H.~Oliveri, T.~Thompson, and A.~Goriely.
    \newblock Front propagation and arrival times in networks with application to neurodegenerative diseases.
    \newblock {\em SIAM Journal on Applied Mathematics}, 83(1):194--224, 2023.
    
    \bibitem{nicholson2000diffusion}
    C.~Nicholson, K.~C. Chen, S.~Hrab{\v{e}}tov{\'a}, and L.~Tao.
    \newblock Diffusion of molecules in brain extracellular space: theory and experiment.
    \newblock {\em Progress in Brain Research}, 125:129--154, 2000.
    
    \bibitem{goodhill1997diffusion}
    G.~J. Goodhill.
    \newblock Diffusion in axon guidance.
    \newblock {\em European Journal of Neuroscience}, 9(7):1414--1421, 1997.
    
    \bibitem{nicholson1998extracellular}
    C.~Nicholson and E.~Sykov{\'a}.
    \newblock Extracellular space structure revealed by diffusion analysis.
    \newblock {\em Trends in Neurosciences}, 21(5):207--215, 1998.
    
    \bibitem{achdou2013qualitative}
    Y.~Achdou, B.~Franchi, N.~Marcello, and M.~C. Tesi.
    \newblock A qualitative model for aggregation and diffusion of $\beta$-amyloid in {A}lzheimer’s disease.
    \newblock {\em Journal of Mathematical Biology}, 67(6-7):1369--1392, 2013.
    
    \bibitem{meisl2016molecular}
    G.~Meisl, J.~B. Kirkegaard, P.~Arosio, T.~C.~T. Michaels, M.~Vendruscolo, C.~M. Dobson, S.~Linse, and T.~P.~J. Knowles.
    \newblock Molecular mechanisms of protein aggregation from global fitting of kinetic models.
    \newblock {\em Nature Protocols}, 11(2):252--272, 2016.
    
    \bibitem{arndt2018structural}
    J.~W. Arndt, F.~Qian, B.~A. Smith, C.~Quan, K.~P. Kilambi, M.~W. Bush, T.~Walz, R.~B. Pepinsky, T.~Bussi{\`e}re, S.~Hamann, et~al.
    \newblock Structural and kinetic basis for the selectivity of aducanumab for aggregated forms of amyloid-$\beta$.
    \newblock {\em Scientific Reports}, 8(1):6412, 2018.
    
    \bibitem{sevigny2016antibody}
    J.~Sevigny, P.~Chiao, T.~Bussi{\`e}re, P.~H. Weinreb, L.~Williams, M.~Maier, R.~Dunstan, S.~Salloway, T.~Chen, Y.~Ling, et~al.
    \newblock The antibody aducanumab reduces {A}$\beta$ plaques in {A}lzheimer’s disease.
    \newblock {\em Nature}, 537(7618):50--56, 2016.
    
    \bibitem{soderberg2023lecanemab}
    L.~S{\"o}derberg, M.~Johannesson, P.~Nygren, H.~Laudon, F.~Eriksson, G.~Osswald, C.~M{\"o}ller, and L.~Lannfelt.
    \newblock Lecanemab, aducanumab, and gantenerumab— binding profiles to different forms of amyloid-beta might explain efficacy and side effects in clinical trials for {A}lzheimer’s disease.
    \newblock {\em Neurotherapeutics}, 20(1):195--206, 2023.
    
    \bibitem{tolar2020aducanumab}
    M.~Tolar, S.~Abushakra, J.~A. Hey, A.~Porsteinsson, and M.~Sabbagh.
    \newblock Aducanumab, gantenerumab, {BAN}2401, and {ALZ}-801 — the first wave of amyloid-targeting drugs for {A}lzheimer’s disease with potential for near term approval.
    \newblock {\em Alzheimer's Research \& Therapy}, 12:1--10, 2020.
    
    \bibitem{withington2022amyloid}
    C.~G. Withington and R.~S. Turner.
    \newblock Amyloid-related imaging abnormalities with anti-amyloid antibodies for the treatment of dementia due to {A}lzheimer's disease.
    \newblock {\em Frontiers in Neurology}, 13:862369, 2022.
    
    \bibitem{cummings2021aducanumab}
    J.~Cummings, P.~Aisen, L.~G. Apostolova, A.~Atri, S.~Salloway, and M.~Weiner.
    \newblock Aducanumab: appropriate use recommendations.
    \newblock {\em The Journal of Prevention of Alzheimer's Disease}, 8:398--410, 2021.
    
    \end{thebibliography}

\end{document}